\documentclass{siamart}

\usepackage{lipsum}
\usepackage{amsfonts}
\usepackage{graphicx}
\usepackage{epstopdf}
\ifpdf
  \DeclareGraphicsExtensions{.eps,.pdf,.png,.jpg}
\else
  \DeclareGraphicsExtensions{.eps}
\fi

\usepackage[normalem]{ulem}
\usepackage{enumitem}
\usepackage{placeins}
\usepackage{amsopn}
\usepackage{amsmath}
\usepackage{upgreek}
\usepackage{color}
\usepackage{empheq}
\usepackage{amsfonts}
\usepackage{tcolorbox}
\usepackage{tikz}
\usepackage{caption}
\usepackage{subcaption}
\usepackage{algorithm}
\usepackage[noend]{algpseudocode}
\usetikzlibrary{shapes,snakes, positioning}
\usetikzlibrary{positioning}

\usepackage{mathtools}
\DeclarePairedDelimiter{\ceil}{\lceil}{\rceil}

\newcommand{\TheAuthors}{N. Alger, P. Chen, and O. Ghattas}

\newcommand{\ShortTitle}{Tensor train from tensor actions}
\newcommand{\LongTitle}{Tensor train construction from tensor actions, with application to compression of large high order derivative tensors}

\newtheorem{defn}{Definition}

\graphicspath{{./}}

\headers{\ShortTitle}{\TheAuthors}

\title{{\LongTitle}\thanks{This research was partially funded by the Department of Energy, Office
		of Science, Office of Advanced Scientific Computing Research,
		Mathematical Multifaceted Integrated Capability Centers (MMICCS)
		program under award DE-SC0019303; the Simons Foundation under award
		560651; the Air Force Office of Scientific Research, Computational
		Mathematics program
		under award FA9550-17-1-0190; and the National Science Foundation, Division of
		Advanced Cyberinfrastructure under award ACI-1550593, and Division of Mathematical Sciences under award DMS-2012453.
}}

\author{
  Nick Alger\thanks{Oden Institute for Computational Engineering and Sciences, The University of Texas at Austin, Austin, TX, USA
  	    (\email{nalger@oden.utexas.edu},
  	    \email{peng@oden.utexas.edu}, \email{omar@oden.utexas.edu})}
  \and
	Peng Chen\footnotemark[2]
  \and
  Omar Ghattas\footnotemark[2]
}

\begin{document}
	
	\maketitle
	
	\begin{abstract}
		We present a method for converting tensors into tensor train format based on actions of the tensor as a vector-valued multilinear function. Existing methods for constructing tensor trains require access to ``array entries'' of the tensor and are therefore inefficient or computationally prohibitive if the tensor is accessible only through its action, especially for high order tensors. Our method permits efficient tensor train compression of large high order derivative tensors for nonlinear mappings that are implicitly defined through the solution of a system of equations. Array entries of these derivative tensors are not directly accessible, but actions of these tensors can be computed efficiently via a procedure that we discuss. Such tensors are often amenable to tensor train compression in theory, but until now no efficient algorithm existed to convert them into tensor train format. We demonstrate our method by compressing a Hilbert tensor of size $41 \times 42 \times 43 \times 44 \times 45$, and by forming high order (up to $5^\text{th}$ order derivatives/$6^\text{th}$ order tensors) Taylor series surrogates of the noise-whitened parameter-to-output map for a stochastic partial differential equation with boundary output.
	\end{abstract}
	

	
	\section{Introduction}
	\label{sec:intro}
	
	To understand this paper, the reader should be familiar with tensors and tensor decompositions \cite{GrasedyckEtAl13,KoldaBader09}, tensor trains \cite{Gelss17,Oseledets11,PerezEtAl06}, and randomized linear algebra \cite{HalkoMartinssonTropp11}.
	Existing tensor train construction methods, including TT-SVD \cite{Oseledets11}, TT-cross \cite{OseledetsTyrtyshnikov10}, Tucker-2/PVD \cite{CichockiLeeOseledetsEtAl16}, modified ALS \cite{HoltzRohwedderSchneider12},  and other related methods \cite{BallaniEtAl13, CoronaRahimianZorin17}, view a $d^\text{th}$ order tensor as a multidimensional array, 
	\begin{equation*}
		T \in \mathbb{R}^{N_1 \times N_2 \times \dots \times N_d},
	\end{equation*}
	and construct a tensor train representation of $T$ via processes that involve accessing many array entries of $T$. This is inefficient in applications where the tensor is only accessible through its \emph{action} as a vector-valued multilinear function. 
	\begin{defn}[Tensor action]
		Let $T$ be a $d^\text{th}$ order tensor. An action of $T$ is a contraction of $T$ with $d-1$ vectors. 
	\end{defn}
	Tensor actions generalize the concept of the action of a matrix on a vector via matrix-vector multiplication. A matrix action takes one vector as input and returns one vector as output. A tensor action takes $d-1$ vectors as input and returns one vector as output. 
	
	To illustrate by example, suppose that $T \in \mathbb{R}^{N_1 \times N_2 \times N_3}$. Then an action of $T$ is an evaluation of one of the following three functions,
	\begin{align*}
	&F_1:\mathbb{R}^{N_2} \times \mathbb{R}^{N_3} \rightarrow \mathbb{R}^{N_1},\\
	&F_2:\mathbb{R}^{N_1} \times \mathbb{R}^{N_3} \rightarrow \mathbb{R}^{N_2},\\
	&F_3:\mathbb{R}^{N_1} \times \mathbb{R}^{N_2} \rightarrow \mathbb{R}^{N_3},
	\end{align*}
	which are defined by
	\begin{align*}
	F_1(v, w)_i :=& \sum_{j=1}^{N_2} \sum_{k=1}^{N_3} T_{i, j, k} v_j w_k, \quad i=1,\dots,N_1, \\
	F_2(u, w)_j :=& \sum_{i=1}^{N_1} \sum_{k=1}^{N_3} T_{i, j, k} u_i w_k, \quad j=1,\dots,N_2,\\ 
	F_3(u, v)_k :=& \sum_{i=1}^{N_1} \sum_{j=1}^{N_2} T_{i, j, k} u_i v_j, \quad k=1,\dots,N_3.
	\end{align*}
	For a $d^\text{th}$ order tensor, there are $d$ functions $F_i$ defined analogously. When we say that a tensor $T$ is only accessible via its actions, we mean that we have algorithms or black-box computer codes that can evaluate the functions $F_i$ on arbitrary input vectors, but we do not have any other information about the entries of $T$. Tensors that are only accessible through their action arise as tensors representing higher order derivatives of inverse and
	optimal control problem objective functions with respect to parameters. More generally, such tensors arise as tensor networks with a tree topology where at least one of the nodes in the network is a matrix inverse.

	
	 If a tensor is only available via actions, it is possible to compress it to tensor train format using existing algorithms such as TT-cross, but this process wastes information.
	 Continuing the example where $T$ is a $3$-tensor, when the existing algorithm wishes to access tensor entry $T_{i,j,k}$, one computes the \emph{fiber} $F_1(e_j, e_k) \in \mathbb{R}^{N_1}$ via a tensor action, then extracts $T_{i,j,k} = F_1(e_j, e_k)_i$ as the $i^\text{th}$ entry of the fiber. Here and throughout the paper we write 
	 \begin{equation*}
	 e_m:=\begin{bmatrix}0 & \dots & 0 & 1 & 0 & \dots & 0\end{bmatrix}^T
	 \end{equation*}
	 to denote a vector of the appropriate length (length $N_2$ and $N_3$ here) with $m^\text{th}$ entry equal to one and all other entries equal to zero. The fiber $F_1(e_j, e_k)$ can be stored and used later when other tensor entries in the fiber are needed (entries $T_{i,j,k}$ with the same $i$ and $j$, but different $k$). However, existing tensor train compression methods access tensor entries in a \emph{scattered} pattern: only a small number of entries in a given fiber are used. A large number of fibers must be computed, each of which requires a tensor action, and most of the entries in the computed fibers are ignored.

	We present an efficient method to construct tensor train representations of tensors that are accessible only through their action. Our method (Section \ref{sec:method}) is based on a randomized algorithm, which has been used in \cite{CheWei19,HuberSchneiderWolf17} for tensor train decomposition. While \cite{CheWei19,HuberSchneiderWolf17} require the ability to perform array operations with the tensor, our method uses only the tensor action. Our key innovation is a method for implicitly computing the action of the remainder of a partially constructed tensor train, even though we do not have access to the remainder directly. This allows us to use a ``peeling process'' to compute successive cores in the tensor train: we use a randomized range finder to compute the first core in the tensor train, then use information from the first core and a randomized range finder to compute the second core, then use information from the first and second cores and a randomized range finder to compute the third core, and so on. This ``peeling process'' continues until we have computed all cores in the tensor train. We were inspired by \cite{LinLuYing11}, in which randomized linear algebra and a peeling process are used to construct hierarchical matrices using only matrix-vector products. With our method, constructing the tensor train requires
	\begin{equation}
	\label{eq:actions_complexity}
	O(\ceil{r/N}dr^2)
	\end{equation}
	tensor actions, plus $O(dNr^2)$ memory for storage and $O(dNr^3)$ operations for linear algebra overhead. Here $r$ is the tensor train rank (the maximum rank of the cores in the tensor train), $N = \max(N_1, \dots, N_d)$, and $\ceil{x}$ is the smallest integer larger than $x$. The tensor actions are trivially parallelizable within each stage of the method. Once the tensor train has been constructed, the tensor may be manipulated efficiently using fast methods for tensor trains. In the case where $r < N$, we have $\ceil{r/N}=1$, so the required number of tensor actions reduces to $O(dr^2)$. Tensor train compression is most effective when the tensor is large (large $N$ and large $d$), and the tensor train rank is small (small $r$), so $r < N$ is a common use case. For tensors that arise in connection with integral or differential equations, $N$ is typically the number of degrees of freedom in a discretization of a continuous function, and is therefore on the order of thousands, millions, or more, while $r$ may be on the order of ten to one hundred.
	
	Our work was motivated by a desire to form high order (derivative order $k > 2$) Taylor series surrogate models of a quantity of interest $\mathcal{F}\in \mathbb{R}^{N_q}$ that depends on a parameter $m\in \mathbb{R}^{N_m}$ implicitly through the solution of a large system of nonlinear equations. Typically, computing a single entry of the quantity of interest, $\left(\mathcal{F}(m)\right)_i$, requires essentially the same amount of work as computing the entire vector $\mathcal{F}(m)$, because one must solve the same system of nonlinear equations in either case. Analogously, computing an individual entry of the $(k+1)^\text{th}$ order tensor
	\begin{equation}
	\label{eq:higher_derivative_tensor81}
	\frac{d^k \mathcal{F}}{dm^k} \in \mathbb{R}^{N_m} \times \dots \times \mathbb{R}^{N_m} \times \mathbb{R}^{N_q}
	\end{equation}
	requires essentially the same work as computing the vector output of a tensor action. Tensor compression algorithms that operate by accessing scattered tensor entries (such as TT-cross) are therefore inefficient here. We show how to compute tensor actions for these higher order derivative tensors in Section \ref{sec:high order_derivative_action}. These high order derivative tensors are often amenable to tensor train compression in principle, but until now no efficient algorithm existed for converting them to tensor train format. Quantized tensor train methods \cite{Oseledets09,OseledetsTyrtyshnikovZamarashkin11,Oseledets10} can be highly effective for building surrogate models of the input to output map for the solution of a linear system \cite{DolgovKhoromskijOseledets12,KazeevKhammashNipSchwab14,KazeevReichmannSchwab13,Khoromskij11}, but require access to tensor array entries, and therefore cannot be used to efficiently approximate higher order derivative tensors that arise when the system is nonlinear.
	
	We demonstrate our method numerically in Section \ref{sec:numerical_results}.\footnote{Code for the methods and numerical results in this paper is available at \url{https://github.com/NickAlger/TensorTrainHigherDerivatives}.} First, we compress a Hilbert tensor of size $41 \times 42 \times 43 \times 44 \times 45$. We show that our method constructs nearly optimal tensor train approximations for this Hilbert tensor compared to TT-SVD and TT-cross. Second, we compress high order derivatives (up to $5^\text{th}$ order derivatives, or $6^\text{th}$ order tensors) of the noise-whitened parameter-to-output map for a stochastic partial differential equation (PDE) with boundary output. We show that high order derivatives can be compressed into tensor train format with a low tensor train rank, and that the number of tensor actions needed to compress the high order derivative tensors is independent of the mesh used to discretize the problem. As the mesh is refined ($N_m \rightarrow \infty$), the required number of tensor actions remains roughly the same. We use these compressed derivative tensors to build Taylor series surrogate models for the noise-whitened parameter-to-output map, and find that including high order terms in the Taylor series yields more accurate surrogate models.

	\subsection{Isomorphism between arrays and multilinear functions}
	
	Tensor contraction establishes an isomorphism between multidimensional arrays and multilinear functions. Given a multilinear function, 
	\begin{equation*}
	T:\mathbb{R}^{N_1} \times \mathbb{R}^{N_2} \times \dots \mathbb{R}^{N_d} \rightarrow \mathbb{R},
	\end{equation*}
	we may form an array representation of the function by applying the function to all possible combinations of standard unit basis vectors $e_m$. Given an array, 
	\begin{equation*}
	T \in \mathbb{R}^{N_1} \times \mathbb{R}^{N_2} \times \dots \mathbb{R}^{N_d},
	\end{equation*}
	we may define a corresponding multilinear function that acts on vectors by contracting those vectors against the modes of the array. This associates each multilinear function with a unique array, and each array with a unique multilinear function. We use the word tensor to refer to both multilinear functions, and multidimensional arrays. When a tensor is viewed as a multilinear function, we use parentheses to denote function arguments, as in $T(u,v,w)$. When a tensor is viewed as an array we use subscripts to denote array entries, as in $T_{i,j,k}$.
	
	We illustrate by continuing the example from the previous section where $T$ is a $3$-tensor. We have
		\begin{equation*}
	T(u,v,w) = \sum_{i=1}^{N_1} \sum_{j=1}^{N_2} \sum_{k=1}^{N_3} T_{i,j,k} u_i v_j w_k
	\end{equation*}
	for vectors $u,v,w$, and
	\begin{equation*}
		T_{i,j,k} := T(e_i, e_j, e_k)
	\end{equation*} 
	for indices $i,j,k$.

	The action of a multilinear function on vectors equals the contraction of the associated array with those vectors. We use the the `$\cdot$' symbol to denote incomplete contraction by currying. That is,
	\begin{equation*}
		T(u,v,~\cdot~)
	\end{equation*}
	denotes the Riesz representation of the linear functional $w \mapsto T(u,v,w)$ with respect to the Euclidean inner product. If $x=T(u,v,~\cdot~)$ then
	\begin{equation*}
		x_k = \sum_{i=1}^{N_1} \sum_{j=1}^{N_2} T_{i,j,k} u_i v_j, \quad k=1,\dots,N_3.
	\end{equation*}
	Using this notation, actions of $T$ in the above example may be written as
	\begin{align*}
		F_1(v,w) =& T(~\cdot~, v, w), \\
		F_2(u,w) =& T(u, ~\cdot~, w), \\
		F_3(u,v) =& T(u, v, ~\cdot~).
	\end{align*}

	\subsection{Tensor train from the multilinear function perspective}
	
	From the conventional array perspective, the cores of a tensor train representation of $T$ are $2$ and $3$-dimensional arrays
	\begin{align*}
		C_1 &\in \mathbb{R}^{N_1 \times r_1}, \\ 
		C_{k+1} &\in \mathbb{R}^{r_k \times N_{k+1} \times r_{k+1}}, \quad k=1,2,\dots, d-2, \\
		C_d &\in \mathbb{R}^{r_{d-1} \times N_d},
	\end{align*}
	such that
	\begin{equation*}
		T_{i_1, i_2, \dots, i_d} = \sum_{i_1=1}^{r_1} \sum_{j_2=1}^{r_2} \dots \sum_{j_{d-1}=1}^{r_{d-1}} \left(C_1\right)_{i_1,j_1} \left(C_2\right)_{j_1,i_2,j_2} \dots \left(C_{d-1}\right)_{j_{d-2},i_{d-1},j_{d-1}} \left(C_{d}\right)_{j_{d-1},i_d}.
	\end{equation*}
	From the multilinear function perspective, this is a factorization of $T$ into the composition of functions,
	\begin{equation*}
		T(x_1, x_2, \dots, x_d) = C_d(C_{d-1}(\dots C_2(C_1(x_1, ~\cdot~), x_2, ~\cdot~)\dots, x_{d-1}, ~\cdot~), x_d),
	\end{equation*}
	where the output from the last mode of each core is used as an input for the first mode of the next core.
	
	\section{Method}
	\label{sec:method}
	
	We construct the cores of the tensor train one at a time. Without loss of generality we assume $d \ge 5$. If $d < 5$, one may proceed as described for cores $1$ through $d-1$, then skip to the method for computing the last core for core $d$. 
	
	The basic idea of the method is as follows. At each step (other than the first), we start with a factorization of $T$ into a partially constructed tensor train, $T_k$, composed with an unknown multilinear remainder function, $R_k$. We construct vectors that, when input into $T_k$, yield the standard unit basis vectors, $e_j$, as output (recall $e_1 = (1,0,\dots,0)$, $e_2 = (0,1,0,\dots,0)$, and so on). This allows us to apply $R_k$ to arbitrary vectors through an indirect process that involves computing actions of $T$. By using a randomized range finding procedure that involves applying $R_k$ to random vectors, we construct the next core in the tensor train. This process repeats until all cores are computed.

	\subsection{Multilinear randomized range finder}
	\label{sec:randomized_range_finder}
	
	In the randomized singular value decomposition (randomized SVD) \cite{HalkoMartinssonTropp11}, one constructs a basis for the range of a matrix by applying that matrix to random input vectors, then orthogonalizing the resulting output vectors. Here we use a multilinear generalization of this idea to construct a basis for the numerical range of a vector valued multilinear function. Similar multilinear randomized range finders have been used in \cite{CheWei19,HuberSchneiderWolf17}.
	
	Let $F$ be a vector-valued multilinear function that takes $n$ vectors as input and returns one vector as output. Let $r$ be the dimension of the numerical range of $F$ (or the desired dimension of an approximation to this range). First, we compute
	\begin{equation*}
		y^{(i)} = F(\omega_1^{(i)}, \omega_2^{(i)}, \dots, \omega_n^{(i)}), \quad i=1,\dots,r+p,
	\end{equation*}
	where $\omega_1^{(i)}$, $\omega_2^{(i)}$, $\dots$, $\omega_n^{(i)}$ are random vectors of the appropriate sizes with independent normally distributed entries, and $p$ is a small oversampling parameter (we use $p=5$). Then we form an orthonormal basis for the span of the $y^{(i)}$ by computing the thin singular value decomposition (SVD) of the matrix that has $y_i$ as its columns. That is,
	\begin{equation*}
		\begin{bmatrix}y^{(1)} & \dots & y^{(r+p)}\end{bmatrix} = U \Sigma V^T.
	\end{equation*}
	Finally, we set $U_r$ to be the matrix consisting of the first $r$ columns of $U$. The span of these $r$ columns approximates the range of $F$.

	\subsection{First core} 
	\label{sec:first_core}
	
	Define $F_1$ to be the following vector valued multilinear map:
	\begin{equation}
	\label{eq:first_action_4242}
	F_1(x_2, \dots, x_d) := T(~\cdot~, x_2, \dots, x_d).
	\end{equation}
	We use the randomized range finding procedure from Section \ref{sec:randomized_range_finder} to compute an orthonormal basis, $U_r$, for the range of $F_1$. Then we set $C_1:=U_r$.

	\subsection{Second core} 
	
	If the span of the columns of $C_1$ accurately captures the range of $F_1$, then $T$ factors into the composition
	\begin{equation}
	\label{eq:first_decomposition}
	T(x_1,x_2,\dots,x_d) = R_1\left(T_1(x_1), x_2, \dots, x_d\right),
	\end{equation}
	where
	\begin{equation*}
		T_1(x_1) := C_1^T x_1,
	\end{equation*}
	and the ``remainder,'' $R_1:\mathbb{R}^{r_1} \times \mathbb{R}^{N_2} \times \dots \times \mathbb{R}^{N_d} \rightarrow \mathbb{R}$, is an unknown multilinear function. To construct the next core, we seek an orthonormal basis for the range of the multilinear function $F_2$ defined as
	\begin{equation}
	\label{eq:R1_12_action}
	F_2(x_3, \dots, x_d) := R_1(~\cdot~, ~\cdot~, x_3, \dots, x_d),
	\end{equation}
	where the output is vectorized.
	
	Let 
	\begin{equation}
	\label{eq:second_core_eta}
	\eta_j :=  \left(C_1\right)_{:, j}
	\end{equation} 
	be the $j^\text{th}$ column of $C_1$ (we use the colon subscript, `$:$', to denote all array entries in a given axis). By orthogonality of $C_1$, we have $C_1^T \eta_j = e_j$, which implies
	\begin{equation}
	\label{eq:T1_ej}
	T_1(\eta_j) = e_j.
	\end{equation}
	Combining \eqref{eq:T1_ej} with \eqref{eq:first_decomposition}, we have
	\begin{equation*}
		R_1(e_j, ~\cdot~, x_3, \dots, x_d) = T(\eta_j, ~\cdot~, x_3, \dots, x_d).
	\end{equation*}
	Stacking these vectors yields
	\begin{equation}
	\label{eq:apply_R1}
	R_1(~\cdot~, ~\cdot~, x_3, \dots, x_d) = \begin{bmatrix}
	T(\eta_1, ~\cdot~, x_3, \dots, x_d) \\
	T(\eta_2, ~\cdot~, x_3, \dots, x_d) \\
	\dots \\
	T(\eta_{r_1}, ~\cdot~, x_3, \dots, x_d)
	\end{bmatrix}.
	\end{equation}
	We may therefore construct an orthonormal basis, $U_r$, for the range of $F_2$ by using the randomized range finding procedure described in Section \ref{sec:randomized_range_finder}. Whenever the randomized range finding procedure requires evaluating the function $F_2$, we do so by forming the right hand side of \eqref{eq:apply_R1}. This, in turn, is done by computing $r_1$ actions of $T$. The second core, $C_2$, is the $r_1 \times N_2 \times r_2$ third order tensor formed by reshaping the $(r_1 N_2) \times r_2$ matrix $U_r$. This process for constructing the second core is summarized in Algorithm \ref{alg:second_core}.
	
	\begin{algorithm}
		\caption{Construction of the $2^\text{nd}$ core}
		\label{alg:second_core}
		\begin{algorithmic}[1]
			\Require Core $C_1$.
			\Ensure Core $C_2$.
			\State Form the vectors $\eta_1, \eta_2, \dots, \eta_{r_1}$ according to \eqref{eq:second_core_eta}.
			\State Compute an orthonormal basis, $U_r \in \mathbb{R}^{(r_1 N_2)\times r_2}$, for the range of $F_2$ using the randomized range finder in Section \ref{sec:randomized_range_finder}, in which \eqref{eq:apply_R1} is used to evaluate $F_2$ as needed within the randomized range finding procedure.
			\State Set $C_2$ to be the $r_1 \times N_2 \times r_2$ reshaped version of $U_r$.
		\end{algorithmic}
	\end{algorithm}
	
	\subsection{Third core} 
	\label{sec:third_core}
	Having computed the first and second cores, the tensor train now factors into the composition
	\begin{equation}
	\label{eq:third_core_equality}
	T(x_1, x_2, x_3,x_4, \dots, x_d) = R_2\left(T_2(x_1, x_2), x_3,x_4, \dots, x_d\right)
	\end{equation}
	where 
	\begin{equation*}
		T_2(x_1, x_2) := C_2\left(C_1^T x_1, x_2, ~\cdot~\right),
	\end{equation*}
	and $R_2:\mathbb{R}^{r_2} \times \mathbb{R}^{N_3} \times \dots \times \mathbb{R}^{N_d} \rightarrow \mathbb{R}$ is an unknown multilinear remainder. Here $C_2\left(C_1^T x_1, x_2, ~\cdot~\right)$ denotes the contraction of $C_2$ with $C_1^T x_1$ in the first mode and $x_2$ in the second mode. 
	
	We must now construct an orthonormal basis for the function $F_3$ defined as
	\begin{equation}
	\label{eq:third_remainder}
	F_3(x_4, \dots, x_d) := R_2 \left(~\cdot~, ~\cdot~, x_4, \dots, x_d\right),
	\end{equation}
	where we view the vectorization of the first two modes of $R_2$ as the output and the remaining modes as the inputs. We seek to find a small number, $\tau$, of vectors $\{\xi_{i}\}_{i=1}^\tau$ and ${\{\eta_{i,j}\}_{i=1}^\tau}_{j=1}^{r_2}$ so that
	\begin{equation}
	\label{eq:original_T_ej_2}
	\sum_{i=1}^\tau T_2(\xi_{i}, \eta_{i,j}) = e_j, \qquad j=1,\dots,r_2,
	\end{equation}
	because then \eqref{eq:third_core_equality} implies
	\begin{equation*}
		R_2\left(e_j, ~\cdot~, x_4, \dots, x_d \right) = \sum_{i=1}^\tau T(\xi_{i}, \eta_{i,j}, ~\cdot~, x_4, \dots, x_d),
	\end{equation*}
	which implies
	\begin{equation}
	\label{eq:R2_T_action}
	R_2\left(~\cdot~, ~\cdot~, x_4, \dots, x_d \right) = \sum_{i=1}^\tau \begin{bmatrix}
	T(\xi_{i}, \eta_{i,1}, ~\cdot~, x_4, \dots, x_d) \\
	T(\xi_{i}, \eta_{i,2}, ~\cdot~, x_4, \dots, x_d) \\
	\vdots \\
	T(\xi_{i}, \eta_{i,r_2}, ~\cdot~, x_4, \dots, x_d)
	\end{bmatrix}.
	\end{equation}
	
	Given $\xi_i$, $\eta_{i,j}$ satisfying \eqref{eq:original_T_ej_2}, we compute an $(r_2 N_3) \times r_3$ orthonormal basis, $U_r$, for the range of $F_3$ using the randomized range finder from Section \ref{sec:randomized_range_finder}. Then we set $C_3$ to be the $r_2 \times N_3 \times r_3$ reshaping of $U_r$ into a $3$-tensor. Whenever the randomized range finder requires evaluating $F_3$, we perform the evaluation by computing the right hand side of \eqref{eq:R2_T_action} for $j=1,\dots,r_2$ by computing $\tau r_2$ actions of $T$.
	
	We now describe how to find $\xi_i$, $\eta_{i,j}$ satisfying \eqref{eq:original_T_ej_2}. We choose 
	\begin{equation}
	\label{eq:psi_ith_column}
	\xi_i := \left(C_1\right)_{:, i}
	\end{equation}
	as the $i^\text{th}$ column of $C_1$. Other choices are possible. We choose \eqref{eq:psi_ith_column} since (a) the randomized range finding procedure is likely to be most accurate for vectors in the span of the first columns of $C_1$, and less accurate for vectors in the span of the later columns, and (b) with the core $C_1$ already computed, the following cores need only be accurate for vectors $x_1$ in the column space of $C_1$.
	
	Given these $\xi_i$ vectors, we may now solve a least squares problem to construct acceptable $\eta_{i,j}$. Let $A_i$, $i=1,\dots,\tau$, be the $N_2 \times r_2$ matrices
	\begin{equation}
	\label{eq:third_core_A_matrix}
	A_i := T_2(\xi_i, ~\cdot~)
	\end{equation}
	formed by contracting the existing partially constructed tensor train with the vectors $\xi_i$. We may write \eqref{eq:original_T_ej_2} as the linear system equation,
	\begin{equation}
	\label{eq:second_linear_system_998}
	\begin{bmatrix}A_1^T & A_2^T & \dots A_\tau^T\end{bmatrix}\begin{bmatrix}
	\eta_{1,j} \\ \eta_{2,j} \\ \vdots \\ \eta_{\tau,j}
	\end{bmatrix} =
	e_j, \quad j = 1, \dots, r_2.
	\end{equation}
	We choose $\tau$ sufficiently large so that \eqref{eq:second_linear_system_998} is underdetermined and therefore generically solvable. Then for $j=1,\dots,r_2$ , we find a least squares solution to \eqref{eq:second_linear_system_998} by QR factorization.
	This linear system has $\tau r_2 N_2$ variables satisfying $r_2^2$ equations, and is therefore underdetermined if $\tau \ge \ceil{r_2/N_2}$. We recommend $\tau = \ceil{r_2/N_2} + 1$, and use this in our numerical results. In the case $r_2 < N_2$, this reduces to $\tau=2$. Choosing $\tau = \ceil{r_2/N_2} + 1$ instead of $\tau = \ceil{r_2/N_2}$ improves performance of the method considerably in our numerical examples, while choosing larger $\tau$ does not. This process for constructing the third core is summarized in Algorithm \ref{alg:third_core}. 
	
	A caveat here is that we cannot choose $\tau$ larger than $r_1$ because there are only $r_1$ fibers $\xi_i$ to choose from $C_1$. If $\tau > r_1$ is required, the algorithm should backtrack and increase $r_1$. Backtracking may be avoided by setting a minimum rank for all cores. In our numerical results, we have $r_i<N_i$ and $\tau = \ceil{r_i/N_i} + 1 = 2$ for all cores, so we ensure that backtracking never occurs by setting a minimum rank of $r_i=2$ for all of the cores.
	
	\begin{algorithm}
		\caption{Construction of the $3^\text{rd}$ core}
		\label{alg:third_core}
		\begin{algorithmic}[1]
			\Require Cores $C_1, C_2$.
			\Ensure Core $C_3$.
			\State Form the vectors $\xi_i, \xi_2, \dots \xi_\tau$ according to \eqref{eq:psi_ith_column}.
			\State Form the matrices $A_1, A_2, \dots, A_\tau$ according to \eqref{eq:third_core_A_matrix}.
			\State Find least-squares solutions to \eqref{eq:second_linear_system_998} for $j=1,\dots,r_2$ to get ${\{\eta_{i,j}\}_{i=1}^\tau}_{j=1}^{r_2}$.
			\State Compute an orthonormal basis, $U_r \in \mathbb{R}^{(r_2 N_3)\times r_3}$, for the range of $F_3$ using the randomized range finder in Section \ref{sec:randomized_range_finder}, in which \eqref{eq:R2_T_action} is used to evaluate $F_3$ as needed within the randomized range finding procedure.
			\State Set $C_3$ to be the $r_2 \times N_3 \times r_3$ reshaped version of $U_r$.
		\end{algorithmic}
	\end{algorithm}

	\subsection{$4^\text{th}$ through $(d-1)^\text{th}$ cores} 
	\label{sec:third_through_dm1}
	The process for constructing the $4^\text{th}$ through $(d-1)^\text{th}$ cores is similar to the process for constructing the $3^\text{rd}$ core. But now there are more modes that must be saturated with specially chosen vectors.
	
	Suppose that we have already computed cores $C_1, C_2, \dots, C_{k}$ for some $k$ in the range $3 \le k \le d-2$, and let $T_k:\mathbb{R}^{N_1 \times N_2 \times \dots \times N_k} \rightarrow \mathbb{R}^{r_k}$ be defined recursively as
	\begin{equation*}
		T_l(x_1,\dots,x_{k-1},x_k) := \begin{cases}
			C_1^T x_1, & l=1, \\
			C_l(T_{l-1}(x_1,\dots,x_{k-1}),x_k, ~\cdot~), & l=2,\dots,k,
		\end{cases}
	\end{equation*}
	where $C_l(T_{l-1}(x_1,\dots,x_{k-1}),x_k, ~\cdot~)$ denotes contraction of $C_l$ with $T_{l-1}(x_1,\dots,x_{k-1})$ in the first mode and $x_k$ in the second mode.
	We have the factorization
	\begin{equation}
	\label{eq:kth_factorization}
	T(x_1,\dots,x_k,x_{k+1},\dots,x_d) = R_k\left(T_k(x_1,\dots,x_k), x_{k+1},\dots,x_d\right),
	\end{equation}
	where $R_k:\mathbb{R}^{r_k} \times \mathbb{R}^{N_{k+1}} \times \dots \times \mathbb{R}^d \rightarrow \mathbb{R}$ is an unknown multilinear function. 
	
	Given vectors $\psi_1, \psi_2, \dots, \psi_{k-2}$, $\{\xi_i\}_{i=1}^\tau$, and ${\{\eta_{i,j}\}_{i=1}^\tau}_{j=1}^{r_{k}}$ satisfying
	\begin{equation}
	\label{eq:general_T_ej_eq}
	\sum_{i=1}^\tau T_k\left(\psi_1, \dots, \psi_{k-2}, \xi_i, \eta_{i,j}\right) = e_j,
	\end{equation}
	we have
	\begin{equation}
	\label{eq:Rk_action_general}
	R_k\left(~\cdot~, ~\cdot~, x_{k+2}, \dots, x_d\right) = \sum_{i=1}^\tau \begin{bmatrix}
	T(\psi_1, \dots, \psi_{k-2}, \xi_i, \eta_{i,1}, ~\cdot~, x_{k+2}, \dots, x_d) \\
	T(\psi_1, \dots, \psi_{k-2}, \xi_i, \eta_{i,2}, ~\cdot~, x_{k+2}, \dots, x_d) \\
	\vdots \\ 
	T(\psi_1, \dots, \psi_{k-2}, \xi_i, \eta_{i,r_k}, ~\cdot~, x_{k+2}, \dots, x_d)
	\end{bmatrix}
	\end{equation}
	by the same argument as the analogous result we presented for the third core.
	We construct an orthonormal basis, $U_r$, for the range of $F_{k+1}$ defined as
	\begin{equation}
	\label{eq:r_generic_74}
	F_{k+1}(x_{k+2}, \dots, x_d) := R_k(~\cdot~,~\cdot~,x_{k+2}, \dots, x_d)
	\end{equation}
	using the randomized range finder from Section \ref{sec:randomized_range_finder}. Within the randomized range finder we use identity \eqref{eq:Rk_action_general} to evaluate $F_{k+1}$ as needed. We set $C_{k+1}$ to be the $r_k \times N_{k+1} \times r_{k+1}$ reshaping of $U_r$. 
	
\begin{figure}
		\centering
		\begin{subfigure}{\textwidth}
			\centering
			\begin{tikzpicture}[scale=0.75, every node/.style={scale=0.75}]
			\node[draw, circle, minimum size=2cm] (t) at (0,0) {$T$};
			\node[draw, circle, minimum size=0.75cm] (tout1) at (0,2) {$x_8$};
			\node[draw, circle, minimum size=0.75cm] (tout2) at (1.414,1.414) {$x_7$};
			\node[draw, circle, minimum size=0.75cm] (tout3) at (2,0) {$x_6$};
			\node[draw, circle, minimum size=0.75cm] (tout4) at (1.414,-1.414) {$x_5$};
			\node[draw, circle, minimum size=0.75cm] (tout5) at (0,-2) {$x_4$};
			\node[draw, circle, minimum size=0.75cm] (tout6) at (-1.414,-1.414) {$x_3$};
			\node[draw, circle, minimum size=0.75cm] (tout7) at (-2,0) {$x_2$};
			\node[draw, circle, minimum size=0.75cm] (tout8) at (-1.414,1.414) {$x_1$};
			
			\draw (t) -- (tout1);
			\draw (t) -- (tout2);
			\draw (t) -- (tout3);
			\draw (t) -- (tout4);
			\draw (t) -- (tout5);
			\draw (t) -- (tout6);
			\draw (t) -- (tout7);
			\draw (t) -- (tout8);
			
			\node (equals) at (3.5,0) {\Large $=$};
			
			\node[draw, circle, minimum size=1.0cm] (C1) at (5,0) {$C_1$};
			\node[draw, circle, minimum size=1.0cm] (C2) at (6.5,0) {$C_2$};
			\node[draw, circle, minimum size=1.0cm] (C3) at (8,0) {$C_3$};
			\node[draw, circle, minimum size=1.5cm] (R3) at (10,0) {$R_3$};
			
			\draw (C1) -- (C2) -- (C3) -- (R3);
			
			\node[draw, circle, minimum size=0.75cm] (u1) at (5,-1.5) {$x_1$};
			\node[draw, circle, minimum size=0.75cm] (u2) at (6.5,-1.5) {$x_2$};
			\node[draw, circle, minimum size=0.75cm] (u3) at (8,-1.5) {$x_3$};
			\node[draw, circle, minimum size=0.75cm] (u4) at (10,-1.5) {$x_4$};
			
			\draw (C1) -- (u1);
			\draw (C2) -- (u2);
			\draw (C3) -- (u3);
			\draw (R3) -- (u4);
			
			\node[draw, circle, minimum size=0.75cm] (u5) at (11.06,-1.06) {$x_5$};
			\node[draw, circle, minimum size=0.75cm] (u6) at (11.5,0) {$x_6$};
			\node[draw, circle, minimum size=0.75cm] (u7) at (11.06,1.06) {$x_7$};
			\node[draw, circle, minimum size=0.75cm] (u8) at (10,1.5) {$x_8$};
			
			\draw (R3) -- (u5);
			\draw (R3) -- (u6);
			\draw (R3) -- (u7);
			\draw (R3) -- (u8);
			
			\end{tikzpicture}
			\caption{Intermediate factorization of $T$ after $3$ cores have been computed.}
			\label{fig:T_equals_CR}
		\end{subfigure}
		\par\bigskip
		\begin{subfigure}{\textwidth}
			\centering
			\begin{tikzpicture}[scale=0.75, every node/.style={scale=0.75}]
			
			\node[draw, circle, minimum size=1.0cm] (C1) at (0,0) {$C_1$};
			\node[draw, circle, minimum size=1.0cm] (C2) at (1.5,0) {$C_2$};
			\node[draw, circle, minimum size=1.0cm] (C3) at (3,0) {$C_3$};
			\node (blank) at (4.5,0) {};
			
			\draw (C1) -- (C2) -- (C3);
			\draw[->] (C3) -- (blank);
			
			\node[draw, circle, minimum size=0.75cm] (u1) at (0,-1.5) {$\psi_1$};
			\node[draw, circle, minimum size=0.75cm] (u2) at (1.5,-1.5) {$\xi_1$};
			\node[draw, circle, minimum size=0.75cm] (u3) at (3,-1.5) {$\eta_j$};
			
			\draw (C1) -- (u1);
			\draw (C2) -- (u2);
			\draw (C3) -- (u3);
			
			\node (equals) at (5.5,0) {\Large $=$};
			
			\node[draw, circle, minimum size=0.75cm, inner sep=0.0] (ej) at (7.0,0) {$e_j$};
			\node (blank2) at (8.4,0) {};
			
			\draw[->] (ej) -- (blank2);
			
			\end{tikzpicture}
			\caption{We solve for $\eta_j$ so that this equality holds.}
			\label{fig:eta_equals_ej}
		\end{subfigure}
		\par\bigskip
		\begin{subfigure}{\textwidth}
			\centering
			\begin{tikzpicture}[scale=0.75, every node/.style={scale=0.75}]
			
			\node[draw, circle, minimum size=2cm] (t) at (0,0) {$T$};
			\node[draw, circle, minimum size=0.75cm, inner sep=0.0] (tout1) at (0,2) {$x_8$};
			\node[draw, circle, minimum size=0.75cm, inner sep=0.0] (tout2) at (1.414,1.414) {$x_7$};
			\node[draw, circle, minimum size=0.75cm, inner sep=0.0] (tout3) at (2,0) {$x_6$};
			\node[draw, circle, minimum size=0.75cm, inner sep=0.0] (tout4) at (1.414,-1.414) {$x_5$};
			\node (tout5) at (0,-2) {};
			\node[draw, circle, minimum size=0.75cm] (tout6) at (-1.414,-1.414) {$\eta_j$};
			\node[draw, circle, minimum size=0.75cm] (tout7) at (-2,0) {$\xi_1$};
			\node[draw, circle, minimum size=0.75cm] (tout8) at (-1.414,1.414) {$\psi_1$};
			
			\draw (t) -- (tout1);
			\draw (t) -- (tout2);
			\draw (t) -- (tout3);
			\draw (t) -- (tout4);
			\draw[->] (t) -- (tout5);
			\draw (t) -- (tout6);
			\draw (t) -- (tout7);
			\draw (t) -- (tout8);
			
			\node (equals) at (3.5,0) {\Large $=$};

			\node[draw, circle, minimum size=0.75cm, inner sep=0.0] (C3) at (5,0) {$e_j$};
			\node[draw, circle, minimum size=1.5cm] (R3) at (6.5,0) {$R_3$};
			
			\draw (C3) -- (R3);
			
			\node (u4) at (6.5,-1.5) {};
			
			\draw[->] (R3) -- (u4);
			
			\node[draw, circle, minimum size=0.75cm, inner sep=0.0] (u5) at (7.56,-1.06) {$x_5$};
			\node[draw, circle, minimum size=0.75cm, inner sep=0.0] (u6) at (8.0,0) {$x_6$};
			\node[draw, circle, minimum size=0.75cm, inner sep=0.0] (u7) at (7.56,1.06) {$x_7$};
			\node[draw, circle, minimum size=0.75cm, inner sep=0.0] (u8) at (6.5,1.5) {$x_8$};
			
			\draw (R3) -- (u5);
			\draw (R3) -- (u6);
			\draw (R3) -- (u7);
			\draw (R3) -- (u8);
			
			\end{tikzpicture}
			\caption{Implicit application of $R_3$ via application of $T$ to specially chosen vectors.}
			\label{fig:Re_equals_T}
		\end{subfigure}
		\caption{Illustration of an intermediate step in the peeling process in the case where $\tau=1$ (typically $\tau>1$, and the left hand sides of (\subref{fig:eta_equals_ej}) and (\subref{fig:Re_equals_T}) are sums of tensors of the form shown here). Given an intermediate tensor train factorization, (\subref{fig:T_equals_CR}), we find certain vectors so that equality holds in (\subref{fig:eta_equals_ej}), then we apply the unknown right factor to vectors implicitly by applying $T$ to certain vectors, (\subref{fig:Re_equals_T}). }
		\label{fig:tensor_diagrams}
	\end{figure}
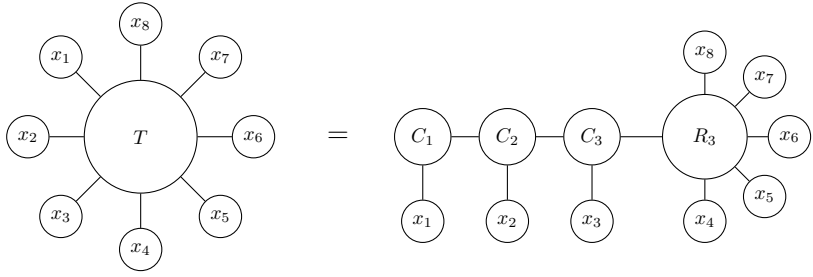
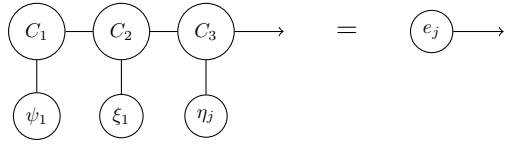
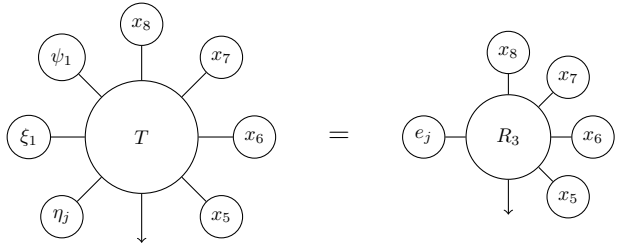

	For the vectors $\psi_1, \psi_2, \dots, \psi_{k-2}$,  we choose $\psi_j$ to be the ``first fibers'' of the corresponding cores $C_j$, which represent the ``highest energy" in each subspace formed by $C_j$, i.e.,
	\begin{equation}
	\label{eq:psi_first_fibers}
	\psi_l := \begin{cases}
	\left(C_1\right)_{:,1}, & l=1, \\
	\left(C_l\right)_{1,:,1}, & l=2, \dots, k-2
	\end{cases}
	\end{equation}
	For the vectors $\{\xi_i\}_{i=1}^\tau$, and ${\{\eta_{i,j}\}_{i=1}^\tau}_{j=1}^{r_{k}}$ satisfying \eqref{eq:general_T_ej_eq}, we use the same process as that for the third core, i.e.,  we specify 
	\begin{equation}
	\label{eq:xi_4th}
	\xi_i := \left(C_{k-1}\right)_{1,:,i}, \quad i=1,\dots,\tau,
	\end{equation}
	and form the $N_k \times r_{k}$ matrices
	\begin{equation}
	\label{eq:A_mtx_def}
	A_i := T_k\left(\psi_1, \dots, \psi_{k-2}, \xi_i, ~\cdot~\right),
	\end{equation}
	and then use a QR factorization to find the least-squares solutions of the linear systems
	\begin{equation}
	\label{eq:generic_least_squares_system_4th}
	\begin{bmatrix}A_1^T & A_2^T & \dots & A_\tau^T\end{bmatrix}\begin{bmatrix}
	\eta_{1,j} \\ \eta_{2,j} \\ \vdots \\ \eta_{\tau,j}
	\end{bmatrix} =
	e_j,
	\end{equation}
	for $j=1,\dots, r_k$. As in the case of the third core, we must have $\tau \ge \ceil{r_{k}/N_k}$ for this system to be solvable; in our numerical experiments we observe good results with $\tau = \ceil{r_{k}/N_k}+1$. This process for computing the $k^\text{th}$ through $(d-1)^\text{th}$ cores is summarized in Algorithm \ref{alg:intermediate_core}. We use graphical tensor notation to illustrate this process in Figure \ref{fig:tensor_diagrams}.
	
	Since we cannot choose $\tau$ larger than $r_{k-1}$, if $\tau > r_{k-1}$ is required, then the algorithm should backtrack and increase $r_{k-1}$. In our numerical results, backtracking is avoided by setting a minimum rank of $r_i=2$ for all of the cores.

	\begin{algorithm}
		\caption{Construction of one of the $4^\text{th}$ through $(d-1)^\text{th}$ cores}
		\label{alg:intermediate_core}
		\begin{algorithmic}[1]
			\Require Cores $C_1, C_2, \dots, C_k$,\quad $3 \le k \le d-2$.
			\Ensure Core $C_{k+1}$.
			\State Form the vectors $\psi_1, \psi_2, \dots \psi_{k-2}$ according to \eqref{eq:psi_first_fibers}.
			\State Form the vectors $\xi_i, \xi_2, \dots \xi_\tau$ according to \eqref{eq:xi_4th}.
			\State Form the matrices $A_1, A_2, \dots, A_\tau$ according to \eqref{eq:A_mtx_def}.
			\State Find least-squares solutions to \eqref{eq:generic_least_squares_system_4th} for $j=1,\dots,r_k$ to get ${\{\eta_{i,j}\}_{i=1}^\tau}_{j=1}^{r_k}$.
			\State Compute an orthonormal basis, $U_r \in \mathbb{R}^{(r_k N_{k+1})\times r_{k+1}}$, for the range of $F_{k+1}$ using the randomized range finder in Section \ref{sec:randomized_range_finder}, in which \eqref{eq:Rk_action_general} is used to evaluate $F_{k+1}$ as needed within the randomized range finding procedure.
			\State Set $C_{k+1}$ to be the $r_k \times N_{k+1} \times r_{k+1}$ reshaped version of $U_r$.
		\end{algorithmic}
	\end{algorithm}

	\subsection{Last core} 
	\label{sec:last_core}
	
	For the last core, we have
	\begin{equation*}
		T(x_1,\dots,x_{d-1},~\cdot~) = R_{d-1}\left(T_{d-1}(x_1,\dots,x_{d-1}), ~\cdot~\right).
	\end{equation*}
	We 	specify the vectors $\psi_1, \dots, \psi_{d-3}$, $\{\xi_i\}_{i=1}^\tau$, and find vectors ${\{\eta_{i,j}\}_{i=1}^\tau}_{j=1}^{r_{d-1}}$ such that
	\begin{equation*}
		\sum_{i=1}^\tau T_{d-1}(\psi_1,\dots,\psi_{d-3}, \xi_i, \eta_{i,j}) = e_j, \quad j = 1, \dots, r_{d-1},
	\end{equation*}
	using the same process as described in Section \ref{sec:third_through_dm1} for previous cores. Then we directly form the last core, $C_d \in \mathbb{R}^{r_{d-1} \times N_d}$, as follows:
	\begin{equation}
	\label{eq:last_pushthrough}
	C_d := \begin{bmatrix}
	R_{d-1}\left(e_1, ~\cdot~\right)^T \\
	R_{d-1}\left(e_2, ~\cdot~\right)^T \\
	\vdots \\
	R_{d-1}\left(e_{r_{d-1}}, ~\cdot~\right)^T
	\end{bmatrix}
	=
	\sum_{i=1}^\tau \begin{bmatrix}
	T(\psi_1,\dots,\psi_{d-3}, \xi_i, \eta_{i,1},~\cdot~)^T \\
	T(\psi_1,\dots,\psi_{d-3}, \xi_i, \eta_{i,2},~\cdot~)^T \\
	\vdots \\
	T(\psi_1,\dots,\psi_{d-3}, \xi_i, \eta_{i,r_{d-1}},~\cdot~)^T
	\end{bmatrix}.
	\end{equation}
	There is no orthogonalization for the last core. This process for constructing the last core is summarized in Algorithm \ref{alg:last_core}.
	
	\begin{algorithm}
		\caption{Construction of the last core}
		\label{alg:last_core}
		\begin{algorithmic}[1]
			\Require Cores $C_1, C_2, \dots, C_{d-1}$.
			\Ensure Core $C_d$.
			\State Form the vectors $\psi_1, \psi_2, \dots \psi_{d-3}$ according to \eqref{eq:psi_first_fibers} with $k=d-1$.
			\State Form the vectors $\xi_i, \xi_2, \dots \xi_\tau$ according to \eqref{eq:xi_4th}, with $k=d-1$.
			\State Form the matrices $A_1, A_2, \dots, A_\tau$ according to \eqref{eq:A_mtx_def}, with $k=d-1$.
			\State Compute least-squares solutions to \eqref{eq:generic_least_squares_system_4th} for $j=1,\dots,r_{d-1}$ to get ${\{\eta_{i,j}\}_{i=1}^\tau}_{j=1}^{r_{d-1}}$.
			\State Construct $C_d$ by evaluating the right hand side of \eqref{eq:last_pushthrough}.
		\end{algorithmic}
	\end{algorithm}

	\subsection{Adaptive range finding}
	\label{sec:adaptive_range_finding}
	
	The randomized range finding procedure may be performed in a sequential manner so that the rank of a core, $r_k$, is found while the orthogonal basis for that core is constructed. Starting with a small $r$ (say, $r=2$), compute $y^{(i)}, i=1,\dots,r+p$ and $U_r$, and compute the error estimate
	\begin{equation}
	\label{eq:error_estimator}
	\mathcal{E} := \max_{i=1,\dots,r+p} \left\|y^{(i)} - U_r U_r^T y^{(i)}\right\|_2.
	\end{equation}
	If $\mathcal{E}$ is less than a predetermined tolerance, set $r_k=r$ and stop the range finding procedure, and reshape $U_r$ to construct the current core. If $\mathcal{E}$ is greater than the tolerance, increase $r$ (say, $r \leftarrow r+1$), compute more vectors $y^{(i)}$, recompute $U_r$, and repeat the process. Error estimator \eqref{eq:error_estimator} generalizes the a posteriori error estimator in Section 4.3 of \cite{HalkoMartinssonTropp11} to tensors. For the numerical results presented in Section \ref{sec:numerical_results}, we fix the rank $r$ beforehand and report results based on other more robust and accurate (but more expensive) methods for computing error.

	\section{Computing the action of high order derivative tensors}
	\label{sec:high order_derivative_action}
	
	In this section we present efficient methods for computing the action of high order derivative tensors for a quantity of interest
	\begin{equation}\label{eq:QoI}
	\mathcal{F}(m,u(m)) \in \mathbb{R}^{N_q},
	\end{equation}
	which depends on a parameter $m \in \mathbb{R}^{N_m}$ implicitly through the solution (state variable) $u \in \mathbb{R}^{N_u}$ of a nonlinear state equation
	\begin{equation}
	\label{eq:nonlinear_equation}
	\mathcal{G}(m,u)=0.
	\end{equation}
	Computation of the entries of the derivative tensor
	\begin{equation}
	\label{eq:higher_derivative_tensor88}
	\frac{d^k \mathcal{F}}{dm^k} \in \mathbb{R}^{N_m} \times \dots \times \mathbb{R}^{N_m} \times \mathbb{R}^{N_q}
	\end{equation}
	for large $N_m$ and $k$ are prohibitive. We show how to compute the action of $S:=\frac{d^k \mathcal{F}}{dm^k}$ as a multilinear function. In Section \ref{sec:output_mode_free}, we show how to compute the action of the derivative tensor when the output mode is free, that is,
	\begin{equation*}
		S(p_1, \dots, p_k, ~\cdot~) = \frac{d^k\mathcal{F}}{dm^k}~p_1~\dots~p_k.
	\end{equation*}
	In Section \ref{sec:derivative_mode_free}, we show how to compute the action of the derivative tensor when a derivative mode is free, that is,
	\begin{equation*}
		S(~\cdot~, p_2, \dots, p_k, q) := q^T \left(\frac{d^k\mathcal{F}}{dm^k}~p_2~p_3~\dots~p_k\right).
	\end{equation*}
	
	Typically one can evaluate directional partial derivatives of $\mathcal{F}$ and $\mathcal{G}$ with respect to $m$ and $u$ with automatic differentiation, but one cannot easily evaluate total derivatives of $\mathcal{F}$ with respect to $m$, because that would require automatically differentiating through the solution procedure for the state equation (e.g., differentiating through an iterative Newton-Krylov solver). The basic idea of this section, therefore, is to convert the problem of computing total derivatives into a problem of computing partial derivatives and solving linear systems.
	We show how one can evaluate the action of $S$ via procedures that involve solving a sequence of linear systems of the form
	\begin{equation}
	\label{eq:incremental_systems}
	0 = \frac{\partial \mathcal{G}}{\partial u}~w + b \qquad \text{and}  \qquad 0 = \left(\frac{\partial \mathcal{G}}{\partial u}\right)^T \!\! w + b,
	\end{equation}
	where constructing the right hand sides $b$ only requires evaluating directional partial derivatives of $\mathcal{F}$ and $\mathcal{G}$, and quantities that have already been computed.  For more on automatic differentiation, we recommend \cite{Naumann12}. A related implementation of computational methods for PDE-constrained high order derivative actions in the finite element package FEniCS \cite{LoggMardalGarth12} can be found in \cite{Maddison19}.

	\subsection{High order derivative actions with the output mode free}
	\label{sec:output_mode_free}
	
	To compute the action of $S$ when the output mode is free, we will repeatedly differentiate $\mathcal{F}$ using the chain rule for total derivatives. This will allow us to express high order total derivatives of $\mathcal{F}$ in terms of partial derivatives of $\mathcal{F}$ and total derivatives of $u$. Repeatedly differentiating the state equation, $0=\mathcal{G}$, in the same manner yields equations that may be solved to determine the required total derivatives of $u$.

	\paragraph{$\text{Zero}^\text{th}$ derivative}
	To compute $\mathcal{F}(m,u(m))$, we solve $\mathcal{G}(m,u)=0$ for $u$, then compute $\mathcal{F}(m,u)$. 
	
	\paragraph{First derivative}
	The chain rule for total derivatives yields
	\begin{equation}
	\label{eq:dqdma}
	\frac{d\mathcal{F}}{dm}~p_a = \frac{\partial \mathcal{F}}{\partial m}~p_a + \frac{\partial \mathcal{F}}{\partial u} \frac{du}{dm}~ p_a = \frac{\partial \mathcal{F}}{\partial m}~p_a +\frac{\partial \mathcal{F}}{\partial u} ~u^{\{a\}},
	\end{equation}
	where $p_a$ is the direction that the derivative is taken in, and
	\begin{equation*}
		u^{\{a\}} := \frac{du}{dm}~ p_a
	\end{equation*}
	is unknown. Differentiating both sides of the state equation, $0 = \mathcal{G}$, yields
	\begin{equation*}
		\frac{d}{dm}\left(0\right) p_a = \frac{d}{dm}\left(\mathcal{G}\right) p_a,
	\end{equation*} 
	or
	\begin{equation}
	\label{eq:first_order_linear_u}
	0 = \frac{\partial \mathcal{G}}{\partial m} ~p_a + \frac{\partial \mathcal{G}}{\partial u} ~u^{\{a\}},
	\end{equation}
	which is a new equation that may be solved for $u^{\{a\}}$. To compute $\frac{d\mathcal{F}}{dm}~p_a$ we solve \eqref{eq:first_order_linear_u} for $u^{\{a\}}$, then evaluate the right hand side of \eqref{eq:dqdma}.
	
	\paragraph{Second derivative}
	Let
	\begin{equation*}
		\mathcal{F}^{\{a\}} := \frac{d\mathcal{F}}{dm}~p_a.
	\end{equation*}
	Using the chain rule for total derivatives again, we have
	\begin{equation}
	\label{eq:second_derivative_u}
	\frac{d^2\mathcal{F}}{dm^2}~p_a ~p_b = \frac{d\mathcal{F}^{\{a\}}}{dm} ~p_b = \frac{\partial \mathcal{F}^{\{a\}}}{\partial m} ~p_b + \frac{\partial \mathcal{F}^{\{a\}}}{\partial u} ~u^{\{b\}} + \frac{\partial \mathcal{F}^{\{a\}}}{\partial u^{\{a\}}} ~u^{\{a,b\}},
	\end{equation}
	where 
	\begin{equation*}
		\qquad u^{\{b\}} := \frac{du}{dm}~ p_b, \quad \text{and} \quad u^{\{a,b\}} := \frac{d^2u}{dm^2}~ p_a ~p_b
	\end{equation*}
	are unknown. We may determine $u^{\{b\}}$ by solving a linear system of the form \eqref{eq:first_order_linear_u}, except with $p_b$ replacing $p_a$. To determine the equation that $u^{\{a,b\}}$ satisfies, we differentiate \eqref{eq:first_order_linear_u} in direction $p_b$. Let 
	\begin{equation*}
		\mathcal{G}^{\{a\}} := \frac{d\mathcal{G}}{dm}~p_a,
	\end{equation*}
	so that \eqref{eq:first_order_linear_u} may be written as $0=\mathcal{G}^{\{a\}}$. We have
	\begin{equation*}
		\frac{d}{dm}\left(0\right) p_b = 	\frac{d}{dm}\left(\mathcal{G}^{\{a\}}\right) p_b
	\end{equation*}
	or
	\begin{equation}
	\label{eq:second_state}
	0 = \frac{\partial \mathcal{G}^{\{a\}}}{\partial m}~p_b + \frac{\partial \mathcal{G}^{\{a\}}}{\partial u}~u^{\{b\}} + \frac{\partial \mathcal{G}^{\{a\}}}{\partial u^{\{a\}}}~u^{\{a,b\}},
	\end{equation}
	which may be solved for $u^{\{a,b\}}$. Once $u^{\{a\}}$ and $u^{\{a,b\}}$ have been determined, we compute $\frac{d^2\mathcal{F}}{dm^2}~p_a ~p_b$ by evaluating the right hand side of \eqref{eq:second_derivative_u}.

	\paragraph{High order derivatives}
	We may repeatedly differentiate $\mathcal{F}$ and $\mathcal{G}$ to construct derivatives of any order. Define
	\begin{align*}
		\mathcal{F}^\alpha :=& \frac{d^k\mathcal{F}}{dm^k}~p_{\alpha_1}~p_{\alpha_2}~\dots~p_{\alpha_k} \\
		\mathcal{G}^\alpha :=& \frac{d^k\mathcal{G}}{dm^k}~p_{\alpha_1}~p_{\alpha_2}~\dots~p_{\alpha_k} \\
		u^\alpha :=& \frac{d^ku}{dm^k}~p_{\alpha_1}~p_{\alpha_2}~\dots~p_{\alpha_k},
	\end{align*}
	where $\alpha$ is a multi-index of $k$ derivative directions (e.g., $\alpha'=\{a,b,b\}$, $k=3$). If $\widetilde{\alpha}:=\{\alpha_1, \alpha_2, \dots, \alpha_{k-1}\}$ is a multi-index of $k-1$ derivative directions created by removing one derivative direction from $\alpha$, then we may generate $\mathcal{F}^{\alpha}$ by differentiating $\mathcal{F}^{\widetilde{\alpha}}$ using the chain rule of total derivatives. This yields
	\begin{align}
		\mathcal{F}^{\alpha} = \frac{d \mathcal{F}^{\widetilde{\alpha}}}{dm}~p_{\alpha_{k}} &= \frac{\partial \mathcal{F}^{\widetilde{\alpha}}}{\partial m}~p_{\alpha_{k}} + \sum_{u^\beta} \frac{\partial \mathcal{F}^{\widetilde{\alpha}}}{\partial u^\beta} \frac{du^\beta}{dm} p_{\alpha_{k}} \nonumber\\
		&= \frac{\partial \mathcal{F}^{\widetilde{\alpha}}}{\partial m}~p_{\alpha_{k}} + 
		\sum_{u^\beta} \frac{\partial \mathcal{F}^{\widetilde{\alpha}}}{\partial u^\beta} u^{\beta \cup \{\alpha_{k}\}}, \label{eq:Falpha005}
	\end{align}
	where the sums are taken over all variables $u^\beta$ that $\mathcal{F}^{\widetilde{\alpha}}$ depends on. If we already have a computer code that computes $\mathcal{F}^{\widetilde{\alpha}}$, then we may use automatic differentiation to create a computer code that computes $\mathcal{F}^{\alpha}$, by using automatic differentiation for each partial derivative in the sum in \eqref{eq:Falpha005}. The code for computing any derivative $\mathcal{F}^\alpha$ may be built by repeated application of this process, by differentiating $\mathcal{F}$ to get $\mathcal{F}^{\{\alpha_1\}}$, differentiating again to get $\mathcal{F}^{\{\alpha_1, \alpha_2\}}$, differentiating again to get $\mathcal{F}^{\{\alpha_1, \alpha_2, \alpha_3\}}$, and so on.
	
	A straightforward inductive argument\footnote{For the base case verify that $\mathcal{F}$ depends on $u$, and for the inductive step use \eqref{eq:Falpha005}.} shows that $\mathcal{F}^\alpha$ depends only on $u^\beta$ for all multiset subsets $\beta \subseteq \alpha$.
	For example, if $\alpha=\{a,b,b\}$, then $\mathcal{F}^\alpha$ depends on $u=u^{\{\}}$, $u^{\{a\}}$, $u^{\{b\}}$, $u^{\{a,b\}}$, $u^{\{b,b\}}$, and $u^{\{a,b,b\}}$. Thus \eqref{eq:Falpha005} may be written as
	\begin{equation}
	\label{eq:Falpha006}
	\mathcal{F}^{\alpha} = \frac{\partial \mathcal{F}^{\widetilde{\alpha}}}{\partial m}~p_{\alpha_{k}} + \sum_{\beta \subseteq \widetilde{\alpha}} \frac{\partial \mathcal{F}^{\widetilde{\alpha}}}{\partial u^\beta} u^{\beta \cup \{\alpha_{k}\}},
	\end{equation}
	where the sum is explicitly written over all multiset subsets $\beta \subset \widetilde{\alpha}$.
	
	We may form $\mathcal{G}^\alpha$ from $\mathcal{G}^{\widetilde{\alpha}}$ in the same manner as
	\begin{equation}
	\label{eq:differentiation_of_G}
	\mathcal{G}^{\alpha} = \frac{\partial \mathcal{G}^{\widetilde{\alpha}}}{\partial m}~p_{\alpha_{k}} + \sum_{\beta \subseteq \widetilde{\alpha}} \frac{\partial \mathcal{G}^{\widetilde{\alpha}}}{\partial u^\beta} u^{\beta \cup \{\alpha_{k}\}}.
	\end{equation}
	Another straightforward argument by induction\footnote{The base case is shown in \eqref{eq:first_order_linear_u}, and the inductive step follows from the chain rule of total derivatives.} shows that there is only one term in the sum on the right hand side of \eqref{eq:differentiation_of_G} that contains $u^{\alpha}$, which takes the form $\frac{\partial \mathcal{G}}{\partial u}u^{\alpha}$. Thus the equation $0=\mathcal{G}^\alpha$, for any $|\alpha|\ge 1$ may be written in the form
	\begin{equation}
	\label{eq:higher_forward99}
	0 = \frac{\partial \mathcal{G}}{\partial u}~u^\alpha + b^\alpha,
	\end{equation}
	where $b^\alpha$ depends on $u^\beta$ for all \emph{strict} multiset subsets $\beta \subset \alpha$, that is, all multiset subsets of $\alpha$, excluding $\alpha$ itself. We may construct $b^\alpha$ by repeatedly automatically differentiating $\mathcal{G}$ via \eqref{eq:differentiation_of_G}, then excluding the term that contains $u^\alpha$ in the result. 
	
\begin{figure}
		\centering
		\begin{subfigure}{0.25\textwidth}
			\centering
			\begin{tikzpicture}
			\node (112) at (0,0) {$\{a,a,a\}$};
			\node (111) at (0,-1) {$\{a,a\}$};
			\node (11) at (0,-2) {$\{a\}$};
			\node (1) at (0,-3) {$\{\}$};
			
			\draw[->] (1) -- (11);
			\draw[->] (11) -- (111);
			\draw[->] (111) -- (112);
			\end{tikzpicture}
			\caption{Symmetric}
		\end{subfigure}
		\begin{subfigure}{0.3\textwidth}
			\centering
			\begin{tikzpicture}
			\node (112) at (1,0) {$\{a,b,b\}$};
			
			\node (111) at (0,-1) {$\{a,b\}$};
			\node (12) at (2,-1) {$\{b,b\}$};
			
			\node (11) at (0,-2) {$\{a\}$};
			\node (2) at (2,-2) {$\{b\}$};
			
			\node (1) at (1,-3) {$\{\}$};

			\draw[->] (1) -- (11);
			\draw[->] (1) -- (2);
			\draw[->] (11) -- (111);
			\draw[->] (2) -- (12);
			\draw[->] (2) -- (111);
			\draw[->] (111) -- (112);
			\draw[->] (12) -- (112);
			\end{tikzpicture}
			\caption{Partially symmetric}
		\end{subfigure}
		\begin{subfigure}{0.35\textwidth}
			\centering
			\begin{tikzpicture}
			\node (abc) at (1,0) {$\{a,b,c\}$};
			
			\node (ab) at (-0.5,-1) {$\{a,b\}$};
			\node (ac) at (1,-1) {$\{a,c\}$};
			\node (bc) at (2.5,-1) {$\{b,c\}$};
			
			\node (a) at (-0.5,-2) {$\{a\}$};
			\node (b) at (1,-2) {$\{b\}$};
			\node (c) at (2.5,-2) {$\{c\}$};
			
			\node (e) at (1,-3) {$\{\}$};

			\draw[->] (e) -> (a);
			\draw[->] (e) -> (b);
			\draw[->] (e) -> (c);
			\draw[->] (a) -> (ab);
			\draw[->] (a) -> (ac);		
			\draw[->] (b) -> (ab);
			\draw[->] (b) -> (bc);
			\draw[->] (c) -> (ac);
			\draw[->] (c) -> (bc);
			\draw[->] (ab) -> (abc);
			\draw[->] (ac) -> (abc);
			\draw[->] (bc) -> (abc);
			\end{tikzpicture}
			\caption{Non-symmetric}
		\end{subfigure}
		\caption{High order quantities of interest, $\mathcal{F}^\alpha$, high order state variables, $u^\alpha$, and high order state equations, $0=\mathcal{G}^\alpha$, are indexed by multi-indices, $\alpha$, consisting of derivative directions. These multi-indices form a lattice ordered by multiset inclusion. A variable $u^\alpha$ depends on the variables $u^\beta$ for all $\beta$ that precede $\alpha$. We compute high order derivatives by working up the lattice, computing $u^\alpha$ at each node by solving $0=\mathcal{G}^\alpha$. The more symmetric the multi-index $\alpha$, the fewer the variables that must be computed. The equations associated with all nodes in a given level of the lattice may be solved in parallel.}
		\label{fig:multiset_lattice}
	\end{figure}
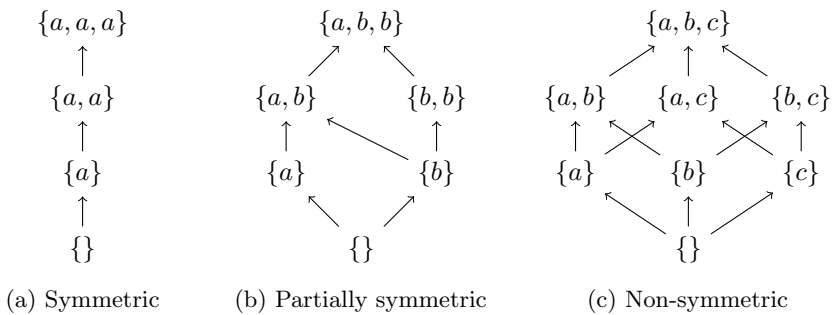

	To determine $u^\alpha$, we need to solve equation \eqref{eq:higher_forward99}, which requires $u^\beta$ for all strict multiset subsets $\beta \subset \alpha$. The dependency structure for the variables $u^\alpha$ and equations $0=\mathcal{G}^\alpha$ therefore forms a bounded lattice consisting of the multiset subsets of $\alpha$, partially ordered by multiset inclusion (see Figure \ref{fig:multiset_lattice}). Thus we may compute $u^\alpha$ by solving the equations $0=\mathcal{G}^\beta$ for the variables $u^\beta$, for all multiset subsets $\beta \subseteq \alpha$, in an order such that each variable $u^\beta$ is solved for after all of the variables corresponding to multiset subsets of $\beta$ have been solved for (e.g., in the order determined by topologically sorting the lattice). This is shown in Algorithm \ref{alg:compute_derivative}. Once these variables are computed, we may compute the desired result,
	\begin{equation}
	S(p_{\alpha_1}, \dots, p_{\alpha_k}, ~\cdot~) = \mathcal{F}^\alpha,
	\end{equation}
	by evaluating the code for $\mathcal{F}^\alpha$ that we generated by repeated automatic differentiation.
	
	\begin{algorithm}
		\caption{Computation of $\mathcal{F}^\alpha = \frac{d^k \mathcal{F}}{dm^k}~p_{\alpha_1}~\dots~p_{\alpha_k}$}\label{alg:compute_derivative}
		\begin{algorithmic}[1]
			\Procedure{compute\_$\mathcal{F}$\_derivative}{$\alpha$}
			\State compute\_$u$\_derivatives($\alpha$)
			\State Construct $\mathcal{F}^\alpha$
			\EndProcedure
			\vspace{1em}
			\Procedure{compute\_$u$\_derivatives}{$\alpha$}
			\If{$u^\alpha$ has already been computed}
			\State Do nothing
			\ElsIf{$\alpha = \{\}$}
			\State Solve $0=\mathcal{G}(m,u)$ for $u$
			\Else
			\For{all multiset subsets $\beta \subset \alpha$ with $|\beta|=|\alpha|-1$}
			\State compute\_$u$\_derivatives($\beta$)
			\EndFor
			\State Construct $b^\alpha$
			\State Solve \eqref{eq:higher_forward99} for $u^\alpha$
			\EndIf
			\EndProcedure
		\end{algorithmic}
	\end{algorithm}
	
	\subsection{High order derivative actions with a derivative mode free}
	\label{sec:derivative_mode_free}
	We may use a similar strategy to compute derivative tensor actions where one of the derivative modes is free. To that end, we define the Lagrangian
	\begin{equation*}
		\mathcal{L}(m,u,\lambda) := q^T \mathcal{F}(u) + \lambda^T\mathcal{G}(m,u),
	\end{equation*}
	where $\lambda \in \mathbb{R}^{N_u}$ is the adjoint variable for enforcing the state equation constraint, and $q\in \mathbb{R}^{N_q}$ is an arbitrary direction that we measure the quantity of interest in. From the theory of Lagrange multipliers, we have the following formula for the first derivative:
	\begin{equation}
	\label{eq:lagrangian_qoi}
	q^T \left(\frac{d\mathcal{F}}{dm}\right) = \frac{\partial \mathcal{L}}{\partial m} = \lambda^T \! \left(\frac{\partial \mathcal{G}}{\partial m}\right),
	\end{equation}
	where $\frac{\partial \mathcal{L}}{\partial m}$ is evaluated at the state $u$, which solves the state equation,
	\begin{equation}
	\label{eq:state}
	0 = \frac{\partial \mathcal{L}}{\partial \lambda} = \mathcal{G}(m,u),
	\end{equation}
	and at the adjoint $\lambda$, which solves the adjoint equation
	\begin{equation}
	\label{eq:adjoint}
	0 = \frac{\partial \mathcal{L}}{\partial u} = q^T \frac{\partial \mathcal{F}}{\partial u} + \lambda^T \! \left(\frac{\partial \mathcal{G}}{\partial u}\right).
	\end{equation}
	This process of computing $q^T \left(\frac{d\mathcal{F}}{dm}\right)$ has the same structure as the process for computing $\mathcal{F}$: we solve a linear system, then evaluate a vector-valued function that depends on the result. Here $\frac{\partial \mathcal{L}}{\partial m}$ takes the place of $\mathcal{F}$, the combined vector $(u,\lambda)$ takes the place of $u$, and the combined state and adjoint system, \eqref{eq:state} and \eqref{eq:adjoint}, takes the place of the state equation. We may therefore compute the desired high order derivatives 
	\begin{equation*}
		q^T \left(\frac{d^k\mathcal{F}}{dm^k}~p_2~p_3~\dots~p_k\right)
	\end{equation*}
	by differentiating \eqref{eq:lagrangian_qoi}, \eqref{eq:state}, and \eqref{eq:adjoint} and solving linear systems repeatedly, the same as was done in Section \ref{sec:output_mode_free}, except with the replacements
	\begingroup
	\renewcommand*{\arraystretch}{1.5}
	\begin{equation*}
		\mathcal{F} \longrightarrow \frac{\partial \mathcal{L}}{\partial m}, \qquad
		u \longrightarrow \begin{bmatrix}
			u \\ \lambda
		\end{bmatrix}, \qquad 
		0=\mathcal{G}(m,u) \longrightarrow \begin{bmatrix}0 \\ 0\end{bmatrix} = \begin{bmatrix}\frac{\partial \mathcal{L}}{\partial \lambda}(m,u) \\  \frac{\partial \mathcal{L}}{\partial u}(m,u,\lambda)\end{bmatrix}.
	\end{equation*}
	\endgroup
	The resulting high order forward equations, 
	\begin{equation*}
		0 = \frac{d^k}{dm^k}\left(\frac{\partial \mathcal{L}}{\partial \lambda}\right)~p_{\alpha_1}~\dots~p_{\alpha_{k-1}},
	\end{equation*}
	are the same as the high order forward equations \eqref{eq:higher_forward99}. The resulting high order adjoint equations,
	\begin{equation*}
		0 = \frac{d^k}{dm^k}\left(\frac{\partial \mathcal{L}}{\partial u}\right)~p_{\alpha_1}~\dots~p_{\alpha_{k-1}},
	\end{equation*}
	take the form
	\begin{equation*}
		0 = \left(\frac{\partial \mathcal{G}}{\partial u}\right)^T \!\! \lambda^\alpha + c^\alpha,
	\end{equation*}
	where $\lambda^\alpha := \frac{d^k \lambda}{dm^k}~p_{\alpha_1}~\dots~p_{\alpha_{k-1}}$, and $c^\alpha$ depends on $u^\beta$ for all multiset subsets $\beta \subseteq \alpha$ and depends on $\lambda^\gamma$ for all strict multiset subsets $\gamma \subset \alpha$. 

	\subsection{Cost to compress derivative tensors}
	
	The method we presented in Section \ref{sec:method} requires $O(\ceil{r/N}dr^2)$ tensor actions to construct a tensor train. But what is the cost of each of these tensor actions? As per the discussion in Section \ref{sec:output_mode_free} and Section \ref{sec:derivative_mode_free}, each high order derivative tensor action requires solving many linear systems of the forms shown in \eqref{eq:incremental_systems}. The number of linear systems that must be solved depends on how symmetric the action of the derivative tensor is. When all the derivative directions are the same, $O(d)$ linear systems must be solved. When all of the derivative directions are distinct, $O(2^d)$ linear systems must be solved. When some directions are the same and some are distinct, an intermediate number of linear systems must be solved. For our method, all derivative directions are distinct, so the number of linear systems that must be solved to construct the tensor train scales as $O(2^d \ceil{r/N} dr^2)$. 
	
	Despite the exponential scaling $2^d$ for a $d$-th order derivative, our method possesses many desirable properties for compressing derivative tensors:
	\begin{itemize}
		\item The $O(2^d)$ linear solves to compute a tensor action are trivially parallelizable to $O(d)$, because the high order variables within each layer of the multiset lattice (e.g., Figure \ref{fig:multiset_lattice}) do not depend on each other, and can therefore be solved for in parallel.
		\item All of the linear systems that must be solved when constructing the tensor train have the same coefficient matrix; they differ only in the right hand side vectors. We therefore construct solvers or preconditioners once, and reuse them for all of the linear systems that must be solved. Constructing solvers or preconditioners is often the most expensive step for solving linear systems.
		\item Up to, say, $5^\text{th}$ or $6^\text{th}$ derivatives (tensor orders $d=6$ or $d=7$), the cost is tractible despite the exponential scaling. This is a substantial improvement over existing methods (e.g., Tucker-based methods), which typically become intractible beyond $2$ or $3$ derivatives.
	\end{itemize}
	
	It is possible to increase the symmetry of the tensor actions within the tensor train construction process, and therefore make the method substantially cheaper, by using the same vector $\psi:=\psi_1=\psi_2=\dots=\psi_{d-3}$, and using the same random vector $\omega_1^{(i)}=\omega_2^{(i)}=\dots=\omega_{d-1}^{(i)}$ for all the derivative modes. These modifications would reduce the number of linear solves per tensor action to $O(d)$, but they may also make the method less robust. We do not use them in our numerical results. At present it is unclear which vector $\psi$ should be chosen. The performance and robustness of the method seem sensitive to the choice of the vectors $\psi_l$.

	\section{Numerical results}
	\label{sec:numerical_results}
	
	In Section \ref{sec:hilbert_tensor}, we compress the Hilbert tensor (a standard test case for tensor compression methods). Tensor entries of the Hilbert tensor are easily computable, so we use this test case to compare the accuracy of our method to that of conventional tensor train compression methods.  
	
	In Section \ref{sec:pde_example}, we compress high order derivative tensors of a quantity of interest that depends on a parameter field implicitly through the solution of a partial differential equation. These derivative tensors are huge for fine meshes, and are available only through their action on vectors; they cannot be compressed to tensor train format efficiently using existing methods. The number of tensor actions for approximating the $k^\text{th}$ derivative tensor using our method is $O\left((k+1)r^2\right)$. For TT-cross, the number of tensor actions would be $O\left(N(k+1)r^2\right)$. We do not compare our method to TT-cross for the derivative tensors in Section \ref{sec:pde_example} because the factor of $N$ makes TT-cross prohibitively expensive. For example, one of the tensors we will compress is a $4^\text{th}$ derivative tensor with $r=30$ and a $40 \times 40$ mesh (part of the Taylor series in Figure \ref{fig:taylor_error_histogram}). This requires approximately $10^4$ tensor actions for our method, and would require approximately $10^7$ tensor actions for TT-cross.

	\subsection{Hilbert tensor}
	\label{sec:hilbert_tensor}
	
	\begin{figure}
		\centering
		\includegraphics[scale=0.6]{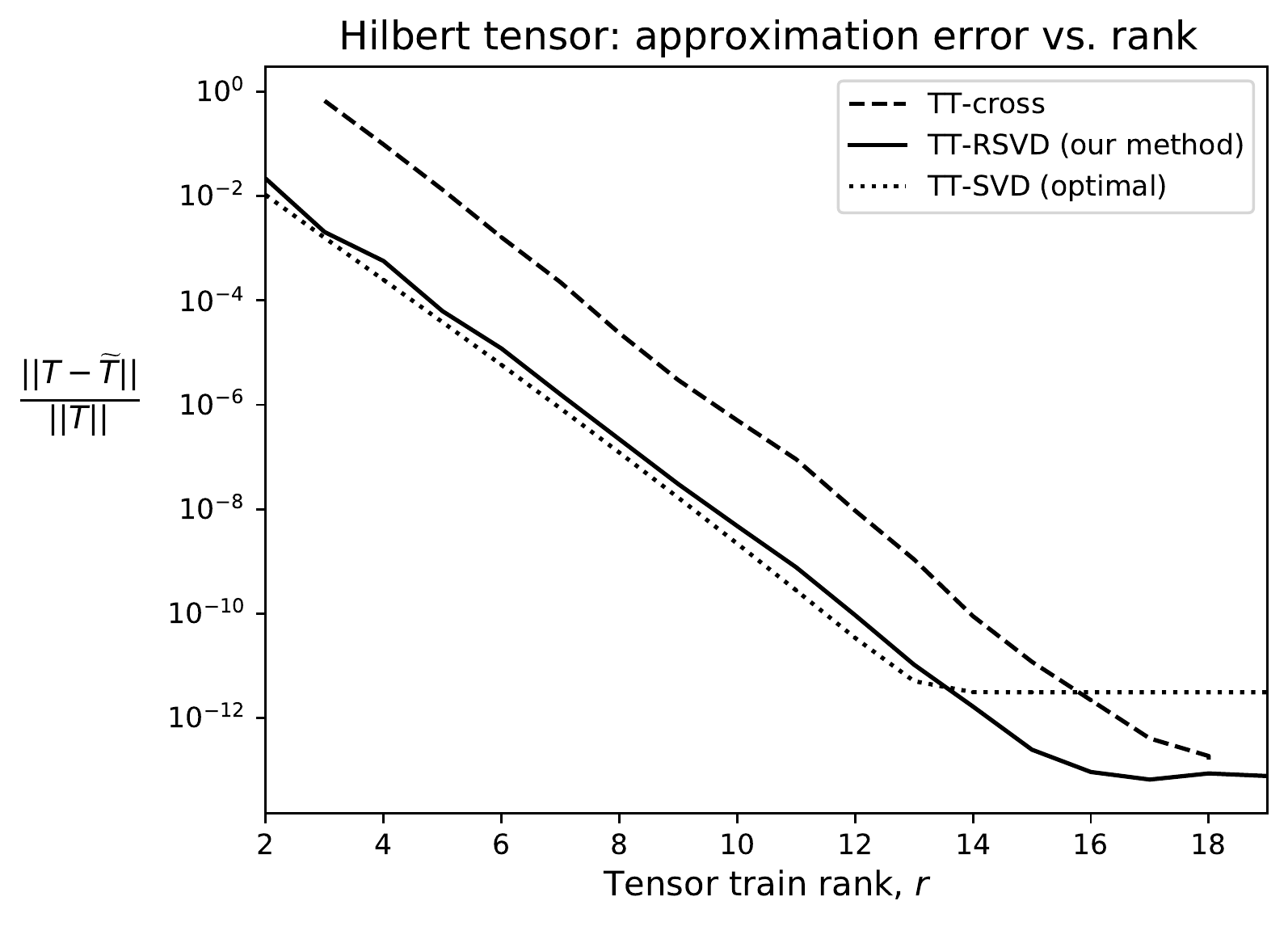}
		\caption{{\bf Hilbert tensor.} Comparison of compression results for the Hilbert tensor using TT-cross, our method (TT-RSVD), and the conventional dense SVD-based algorithm (TT-SVD).}
		\label{fig:hilbert_results}
	\end{figure}
	
	Here we compress the $41 \times 42 \times 43 \times 44 \times 45$ Hilbert tensor, $T$, with entries
	\begin{equation*}
		T\left[i_1, i_2, i_3, i_4, i_5\right] = \frac{1}{i_1 + i_2 + i_3 + i_4 + i_5}.
	\end{equation*}
	In Figure \ref{fig:hilbert_results}, we compare our method with the conventional dense SVD-based method for constructing tensor trains, and with TT-cross \cite{OseledetsTyrtyshnikov10}. For TT-cross, we use the \verb|dmrg_cross()| function in the TT-Toolbox software package\footnote{Version 2.2.2, downloaded from https://github.com/oseledets/TT-Toolbox} \cite{SavostyanovOseledets11}. We compute tensor train approximations, $\widetilde{T}$, of $T$ for a sequence of ranks $r$. The same rank, $r$, is used for all cores. Then we reconstitute full tensors from the tensor train approximations, and compute the relative errors $||T - \widetilde{T}||/||T||$, where $||\cdot||$ is the Frobenius norm (square root of the sum of squares of all tensor entries). We plot the relative error against the rank, $r$. Both TT-cross and our method (TT-RSVD) achieve compression results that are close to the conventional dense SVD-based method (TT-SVD), which is theoretically optimal if one ignores rounding error due to numerical precision. Our method performs slightly better than TT-cross.

	\subsection{High order derivatives of PDE-dependent quantity of interest}
	\label{sec:pde_example}
	
	Here we use our method to build Taylor series surrogate models for the noise-whitened parameter-to-output map for a nonlinear stochastic PDE with boundary output. We define the state equation, $0=\mathcal{G}(m,u)$, to be the following inhomogeneous nonlinear reaction-diffusion equation with Neumann boundary conditions: 
	\begin{equation*}
		\begin{cases}
			-\nabla \cdot e^m \nabla u + u^3 = \rho & \text{in } \Omega, \\
			\nu \cdot e^m \nabla u = 0 & \text{on } \partial \Omega.
		\end{cases}
	\end{equation*}
	The parameter, $m \sim N(0,C)$, is a Gaussian random field with mean zero and covariance
	\begin{equation*}
		C = \left(-\Delta + I\right)^{-2},
	\end{equation*}
	where $\Delta$ is the Laplacian defined in $\Omega$ with Neumann boundary conditions along $\partial \Omega$, and $I$ is the identity operator. The domain is the unit square $\Omega=[0,1]^2$, $\nu$ is the normal to the boundary, and the source term, $\rho: \Omega \to \mathbb{R}$, is given by
	\begin{equation*}
		\rho(c) = \exp\left(\frac{\left\|c - (0.5, 0.5)\right\|^2}{2\left(0.2\right)^2}\right)
	\end{equation*}
	where we use $c \in \Omega$ to denote a spatial point in $\Omega$, to distinguish from $x$ as the vector for the tensor used throughout the paper. We choose the quantity of interest, $\mathcal{F}$, to be the trace of $u$ along the boundary as
	\begin{equation*}
		q^T \mathcal{F}(u) := \int_{\partial \Omega} u~ q ~ds.
	\end{equation*}
	We assume $x \sim N(0, I)$ is a spatial white noise, or a Gaussian distribution with mean zero and covariance $I$, so that
	\begin{equation*}
		p := C^{1/2} x = \left(-\Delta + I\right)^{-1} x
	\end{equation*}
	has a Gaussian distribution with mean zero and covariance $C$. The problem is discretized with $P^1$ finite elements on a mesh of triangles arranged in a regular grid. 
	
	We use our method to approximate the noise-whitened $k^\text{th}$ derivative tensor $T$, defined as
	\begin{equation*}
		T(x_1, \dots, x_k, q) := S(p_1, \dots, p_k, q) = q^T\left(\frac{d^k\mathcal{F}}{dm^k}~p_1~\dots~p_k\right),
	\end{equation*}
	with a rank-$r$ tensor train, $\widetilde{T}$. Note that for a $k^\text{th}$ order derivative, the tensor order is $d=k+1$. 
	By performing this tensor train compression for several of these derivative tensors, we approximate the noise-whitened parameter-to-output map,
	\begin{equation*}
		f(x) := \mathcal{F}\left(u\left(C^{1/2}x\right)\right),		
	\end{equation*}
	with a tuncated Taylor series
	\begin{equation}
	\label{eq:truncated_taylor_series_345}
	f_k(x) =  f(0) + \frac{df}{dx}(0)x + \frac{1}{2}\frac{d^2 f}{dx^2}(0)x^2 + \dots + \frac{1}{k!}\frac{d^k f}{dx^k}(0)x^k.
	\end{equation}
	We write $\widetilde{f}_k$ to denote the truncated Taylor series $f_k$, with the derivative tensors in $f_k$ replaced by their tensor train approximations. 
	
	Surrogate models based on truncation of Taylor series up to the linear term are common in Bayesian inversion, stochastic optimization, and model reduction, where they can be used directly as approximations of $\mathcal{F}(m)$, or used in Markov chain Monte Carlo proposals for Bayesian inversion as variance reduction devices \cite{BashirWillcoxGhattasEtAl08, BigoniZahmSpantiniEtAl2019greedy, Bui-ThanhGhattas2012scaled,
	Bui-ThanhGhattasMartinEtAl13, CuiLawMarzouk16, IsaacPetraStadlerEtAl15, PetraMartinStadlerEtAl14, SpantiniSolonenCuiEtAl15, ZhuLiFomelEtAl16}. Their advantage is that once constructed, no more PDEs need to be solved during the sampling or optimization process. Methods that truncate the Taylor series after the quadratic term have been investigated in \cite{AlexanderianPetraStadlerEtAl17,ChenGhattas19a,ChenVillaGhattas17,ChenVillaGhattas19,ChenWuChenEtAl19a}. In these papers, the second derivative (a third order tensor) is not compressed and stored; rather the action of this tensor is used in other ways. In \cite{BonizzoniNobile14, BonizzoniNobileKressner14}, high order derivatives for the parameter-to-solution map for the log-normal linear Darcy problem are compressed into tensor train format, and Taylor series are constructed. Their algorithm depends specifically on the linear Darcy problem with the solution map as the quantity of interest. With the methods in this paper, we now have a general-purpose algorithm for compressing and storing high order derivative tensors for nonlinear problems with arbitrary quantity of interest. 
	
	\begin{table}
		\centering
		\begingroup
		\renewcommand*{\arraystretch}{1.0}
		\begin{tabular}{c|ccccc|ccccc|}
			& \multicolumn{5}{|c|}{$\epsilon = 10^{-2}$} & \multicolumn{5}{|c|}{$\epsilon=10^{-3}$} \\
			mesh & $k=1$ & $2$ & $3$ & $4$ & $5$ & $k=1$ & $2$ & $3$ & $4$ & $5$ \\
			\hline
			$10\times 10$  & 9  & 8  & 8  & 8  & 11  & 19  & 19  & 20  & 23  & 28  \\
			$20\times 20$  & 8  & 9  & 9  & 10 & 11  & 23  & 22  & 24  & 28  & 35  \\
			$30\times 30$  & 8  & 8  & 9  & 10 & 10  & 22  & 24  & 25  & 29  & 32  \\
			$40\times 40$  & 9  & 9  & 9  & 8  & 11  & 23  & 25  & 27  & 29  & 34  \\
			$50\times 50$  & 7  & 9  & 9  & 11 & 11  & 26  & 23  & 29  & 28  & 34  \\
			$60\times 60$  & 8  & 8  & 10 & 9  & 11  & 24  & 25  & 27  & 30  & 35  \\
			$70\times 70$  & 9  & 9  & 8  & 8  & 12  & 27  & 23  & 27  & 30  & 35  \\
			$80\times 80$  & 9  & 9  & 8  & 11 & 11  & 26  & 25  & 27  & 32  & 38  \\
		\end{tabular}
		\endgroup
		\caption{\textbf{Mesh scalability.} Rank $r$ required to compress derivative tensors to relative error tolerance $\sigma_1(T - \widetilde{T}) / \sigma_1(T) < \epsilon$, where $\epsilon=10^{-2}$ or $\epsilon=10^{-3}$, for a variety of mesh sizes and derivative orders. Meshes are triangles arranged in a regular rectilinear grid, ranging from $10 \times 10$ to $80 \times 80$. Derivative orders range from the $1^\text{st}$ derivative (the Jacobian, a $2^\text{nd}$ order tensor) to $5^\text{th}$ derivative (a $6^\text{th}$ order tensor). For derivative orders $k=2$ through $k=5$, $r$ is the tensor-train rank required using our method. For $k=1$ (Jacobian), $r$ is the matrix rank required using randomized SVD.}  
		\label{tbl:mesh_scalability}
	\end{table}
	
	In Table \ref{tbl:mesh_scalability}, we report the rank required to achieve the relative error 
	\begin{equation*}
		\frac{\sigma_1\left(T - \widetilde{T}\right)}{\sigma_1\left(T\right)} < \epsilon,
	\end{equation*}
	where $\epsilon=10^{-2}$ or $\epsilon=10^{-3}$, and
	\begin{equation}
	\label{eq:sigma_1_def}
	\sigma_1(T) := \max_{\|x\|=1} \|T(x,\dots,x, ~\cdot~)\|.
	\end{equation}
	We use $\sigma_1$ to measure the error because the form of the argument being maximized in \eqref{eq:sigma_1_def} mimics the way the tensor is used within the Taylor series, \eqref{eq:truncated_taylor_series_345}. 
	Computing $\sigma_1$ is NP-hard \cite{HillarLim13}, so we estimate $\sigma_1$ by applying a shifted symmetric high order power method \cite{KoldaMayo11} to the function $x \mapsto T(x,\dots,x, T(x, \dots, x, ~\cdot~))$, with $5$ random initial guesses. 
	
	We show results for derivative orders $k=1$ through $k=5$ (tensor orders $d=2$ through $d=6$) and meshes ranging from $10 \times 10$ to $80\times 80$. For $k=2$ through $k=5$, the reported rank $r$ is the tensor-train rank required using our method. For $k=1$, $r$ is the matrix rank required using randomized SVD. Our results show that the tensor train approximation is mesh-scalable for this problem. As the mesh is refined, the required tensor train rank remains roughly constant. 
	
	\begin{figure}
		\begin{subfigure}[b]{0.55\textwidth}
			\centering
			\includegraphics[scale=0.5]{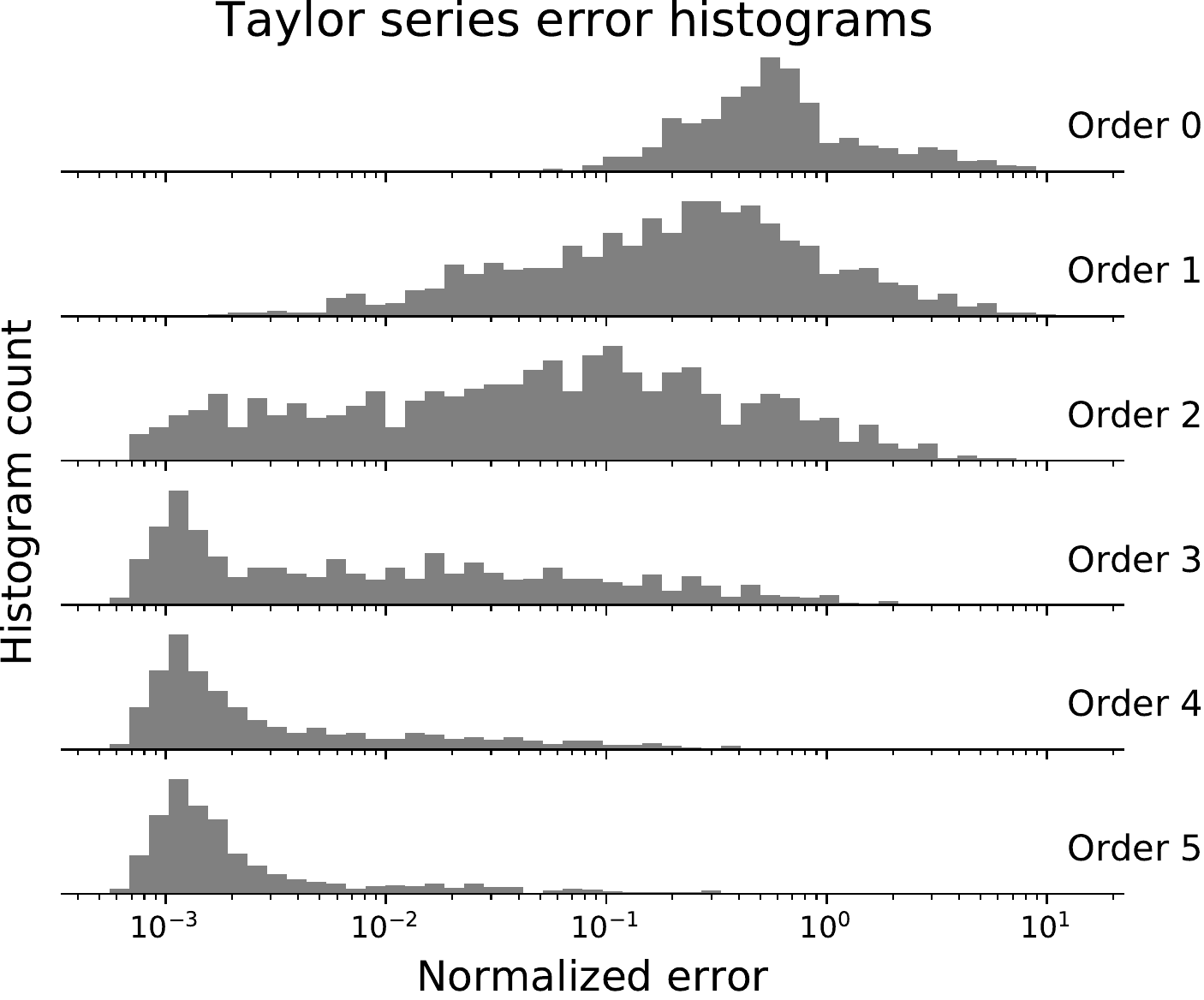}
			\caption{}
		\end{subfigure}
		~
		\begin{subfigure}[b]{0.4\textwidth}
			\centering
			\begin{tabular}{c|c|c}
				order & $\substack{\text{error} \\ \text{mean}}$ & $\substack{\text{error} \\ \text{standard} \\ \text{deviation}}$ \\
				\hline
				0 & 1.0 & 1.4 \\
				1 & 0.55 & 0.98 \\
				2 & 0.24 & 0.54 \\
				3 & 0.078 & 0.21 \\
				4 & 0.019 & 0.069 \\
				5 & 0.014 & 0.075
			\end{tabular}
			\caption{}
		\end{subfigure}
		\caption{{\bf Taylor series error.} (a) Error histograms for the normalized error $\|f(x) - \widetilde{f}_k(x)\| / \mathbb{E}\left(\|f - f(0)\|\right)$, for $1000$ random samples $x \sim N(0,I)$, for Taylor series' of order $0$ through $5$, using a $40 \times 40$ mesh. Rank $30$ is used for all derivative tensors. The high order derivative tensors ($2^\text{nd}$ through $5^\text{th}$ derivatives) are approximated with our method. The Jacobian ($1^\text{st}$ derivative) is approximated with randomized SVD. (b) The means and standard deviations of these normalized error distributions.}
		\label{fig:taylor_error_histogram}
	\end{figure}
	
	In Figure \ref{fig:taylor_error_histogram} we show histograms for the normalized error
	\begin{equation*}
		\frac{\|f(x) - \widetilde{f}_k(x)\|}{\mathbb{E}\left(\|f - f(0)\|\right)},
	\end{equation*}
	for $1000$ random samples $x \sim N(0,I)$ (i.e., the normalized error in the Taylor approximation of the noise-whitened parameter-to-output map, for samples drawn from the noise-whitened parameter distribution). We show Taylor series of order $0$ through $5$. Rank $30$ approximation is used for all derivative tensors. We use a $40 \times 40$ mesh. The high order derivative tensors ($2^\text{nd}$ through $5^\text{th}$ derivative tensors) are approximated with our method, while the Jacobian ($1^\text{st}$ derivative) is approximated with randomized SVD. Constructing the $5^\text{th}$ order Taylor series takes approximately 25 minutes on a laptop (HP model 17-CA1031DX).
	
	Including high order derivatives increases the approximation accuracy. As the order of the Taylor series increases, the error in the approximation decreases. Also, as the order of the Taylor series increases, the normalized error distribution concentrates near $10^{-3}$, which is roughly the error due to the rank-$30$ tensor train approximations of the derivative tensors.

	\section{Conclusion}
	We developed a new randomized algorithm for tensor train decomposition of tensors whose entries cannot be directly accessed or are very expensive to compute. Our method requires only the tensor's action on given vectors and uses a ``peeling process" to successively compute the tensor train cores. This process requires $O(\ceil{r/N}dr^2)$ tensor actions for a $d$-th order tensor with maximum rank $r$ of the cores. We demonstrated the accuracy of this method compared to the conventional TT-SVD and TT-cross methods for a Hilbert tensor. Moreover, we applied this method to construct tensor train decompositions of PDE-constrained high order derivative tensors, for which we derived an efficient scheme to compute the action of arbitrary order derivative tensors. Our method now enables one to use Taylor series truncated to high order terms for uncertainty quantification \cite{BonizzoniNobile14, BonizzoniNobileKressner14}, Bayesian inversion \cite{ChenVillaGhattas17, ChenWuChenEtAl19a}, stochastic optimization \cite{AlexanderianPetraStadlerEtAl17, ChenVillaGhattas19}, and model reduction \cite{BachCegliaSongEtAl19, ChenGhattas19a}. Furthermore, other fields offer promising potential applications of the tensor action based tensor train decomposition, such as neural networks \cite{DengSunQianEtAl19, YangKrompassTresp17} and model constrained high dimensional sampling and integration \cite{DolgovAnaya-IzquierdoFoxEtAl19}.

	\section*{Acknowledgements}
	
	We thank Alen Alexandrian, J.J. Alger, Josh Chen, Tan Bui-Thanh, Andrew Potter, Keyi Wu, and Qiwei Zhan for helpful discussions.
	
	\FloatBarrier
	\bibliographystyle{siam}
	\bibliography{references}

\begin{thebibliography}{10}

\bibitem{AlexanderianPetraStadlerEtAl17}
{\sc A.~Alexanderian, N.~Petra, G.~Stadler, and O.~Ghattas}, {\em Mean-variance
  risk-averse optimal control of systems governed by {PDEs} with random
  parameter fields using quadratic approximations}, SIAM/ASA Journal on
  Uncertainty Quantification, 5 (2017), pp.~1166--1192.
\newblock arXiv preprint arXiv:1602.07592.

\bibitem{BachCegliaSongEtAl19}
{\sc C.~Bach, D.~Ceglia, L.~Song, and F.~Duddeck}, {\em Randomized low-rank
  approximation methods for projection-based model order reduction of large
  nonlinear dynamical problems}, International Journal for Numerical Methods in
  Engineering, 118 (2019), pp.~209--241.

\bibitem{BallaniEtAl13}
{\sc J.~Ballani, L.~Grasedyck, and M.~Kluge}, {\em Black box approximation of
  tensors in hierarchical {T}ucker format}, Linear algebra and its
  applications, 438 (2013), pp.~639--657.

\bibitem{BashirWillcoxGhattasEtAl08}
{\sc O.~Bashir, K.~Willcox, O.~Ghattas, B.~van Bloemen~Waanders, and J.~Hill},
  {\em Hessian-based model reduction for large-scale systems with initial
  condition inputs}, International Journal for Numerical Methods in
  Engineering, 73 (2008), pp.~844--868.

\bibitem{BigoniZahmSpantiniEtAl2019greedy}
{\sc D.~Bigoni, O.~Zahm, A.~Spantini, and Y.~Marzouk}, {\em Greedy inference
  with layers of lazy maps}, arXiv preprint arXiv:1906.00031,  (2019).

\bibitem{BonizzoniNobile14}
{\sc F.~Bonizzoni and F.~Nobile}, {\em Perturbation analysis for the {D}arcy
  problem with log-normal permeability}, SIAM/ASA Journal on Uncertainty
  Quantification, 2 (2014), pp.~223--244.

\bibitem{BonizzoniNobileKressner14}
{\sc F.~Bonizzoni, F.~Nobile, and D.~Kressner}, {\em Tensor train approximation
  of moment equations for the log-normal {D}arcy problem}, tech. rep., Mathicse
  Technical Report, 2014.

\bibitem{Bui-ThanhGhattas2012scaled}
{\sc T.~Bui-Thanh and O.~Ghattas}, {\em A scaled stochastic {N}ewton algorithm
  for {M}arkov chain {M}onte {C}arlo simulations}, SIAM Journal on Uncertainty
  Quantification,  (2012), pp.~1--25.

\bibitem{Bui-ThanhGhattasMartinEtAl13}
{\sc T.~Bui-Thanh, O.~Ghattas, J.~Martin, and G.~Stadler}, {\em A computational
  framework for infinite-dimensional {B}ayesian inverse problems {P}art {I}:
  {T}he linearized case, with application to global seismic inversion}, SIAM
  Journal on Scientific Computing, 35 (2013), pp.~A2494--A2523.

\bibitem{CheWei19}
{\sc M.~Che and Y.~Wei}, {\em Randomized algorithms for the approximations of
  {T}ucker and the tensor train decompositions}, Advances in Computational
  Mathematics, 45 (2019), pp.~395--428.

\bibitem{ChenGhattas19a}
{\sc P.~Chen and O.~Ghattas}, {\em Hessian-based sampling for high-dimensional
  model reduction}, International Journal for Uncertainty Quantification, 9
  (2019).

\bibitem{ChenVillaGhattas17}
{\sc P.~Chen, U.~Villa, and O.~Ghattas}, {\em Hessian-based adaptive sparse
  quadrature for infinite-dimensional {B}ayesian inverse problems}, Computer
  Methods in Applied Mechanics and Engineering, 327 (2017), pp.~147--172.

\bibitem{ChenVillaGhattas19}
{\sc P.~Chen, U.~Villa, and O.~Ghattas}, {\em Taylor approximation and variance
  reduction for {PDE}-constrained optimal control under uncertainty}, Journal
  of Computational Physics, 385 (2019), pp.~163--186.

\bibitem{ChenWuChenEtAl19a}
{\sc P.~Chen, K.~Wu, J.~Chen, T.~O'Leary-Roseberry, and O.~Ghattas}, {\em
  Projected {S}tein variational {N}ewton: {A} fast and scalable {B}ayesian
  inference method in high dimensions}, Advances in Neural Information
  Processing Systems,  (2019).

\bibitem{CichockiLeeOseledetsEtAl16}
{\sc A.~Cichocki, N.~Lee, I.~Oseledets, A.-H. Phan, Q.~Zhao, D.~P. Mandic,
  et~al.}, {\em Tensor networks for dimensionality reduction and large-scale
  optimization: {P}art 1 low-rank tensor decompositions}, Foundations and
  Trends{\textregistered} in Machine Learning, 9 (2016), pp.~249--429.

\bibitem{CoronaRahimianZorin17}
{\sc E.~Corona, A.~Rahimian, and D.~Zorin}, {\em A tensor-train accelerated
  solver for integral equations in complex geometries}, Journal of
  Computational Physics, 334 (2017), pp.~145--169.

\bibitem{CuiLawMarzouk16}
{\sc T.~Cui, K.~Law, and Y.~Marzouk}, {\em Dimension-independent
  likelihood-informed {MCMC}}, Journal of Computational Physics, 304 (2016),
  pp.~109--137.

\bibitem{DengSunQianEtAl19}
{\sc C.~Deng, F.~Sun, X.~Qian, J.~Lin, Z.~Wang, and B.~Yuan}, {\em Tie:
  energy-efficient tensor train-based inference engine for deep neural
  network}, in Proceedings of the 46th International Symposium on Computer
  Architecture, 2019, pp.~264--278.

\bibitem{DolgovAnaya-IzquierdoFoxEtAl19}
{\sc S.~Dolgov, K.~Anaya-Izquierdo, C.~Fox, and R.~Scheichl}, {\em
  Approximation and sampling of multivariate probability distributions in the
  tensor train decomposition}, Statistics and Computing,  (2019), pp.~1--23.

\bibitem{DolgovKhoromskijOseledets12}
{\sc S.~V. Dolgov, B.~N. Khoromskij, and I.~V. Oseledets}, {\em Fast solution
  of parabolic problems in the tensor train/quantized tensor train format with
  initial application to the fokker--planck equation}, SIAM Journal on
  Scientific Computing, 34 (2012), pp.~A3016--A3038.

\bibitem{Gelss17}
{\sc P.~Gel{\ss}}, {\em The Tensor-Train Format and Its Applications: Modeling
  and Analysis of Chemical Reaction Networks, Catalytic Processes, Fluid Flows,
  and Brownian Dynamics}, PhD thesis, 2017.

\bibitem{GrasedyckEtAl13}
{\sc L.~Grasedyck, D.~Kressner, and C.~Tobler}, {\em A literature survey of
  low-rank tensor approximation techniques}, GAMM-Mitteilungen, 36 (2013),
  pp.~53--78.

\bibitem{HalkoMartinssonTropp11}
{\sc N.~Halko, P.~G. Martinsson, and J.~A. Tropp}, {\em Finding structure with
  randomness: {P}robabilistic algorithms for constructing approximate matrix
  decompositions}, SIAM Review, 53 (2011), pp.~217--288.

\bibitem{HillarLim13}
{\sc C.~J. Hillar and L.-H. Lim}, {\em Most tensor problems are {NP}-hard},
  Journal of the ACM (JACM), 60 (2013), pp.~1--39.

\bibitem{HoltzRohwedderSchneider12}
{\sc S.~Holtz, T.~Rohwedder, and R.~Schneider}, {\em The alternating linear
  scheme for tensor optimization in the tensor train format}, SIAM Journal on
  Scientific Computing, 34 (2012), pp.~A683--A713.

\bibitem{HuberSchneiderWolf17}
{\sc B.~Huber, R.~Schneider, and S.~Wolf}, {\em A randomized tensor train
  singular value decomposition}, in Compressed Sensing and its Applications,
  Springer, 2017, pp.~261--290.

\bibitem{IsaacPetraStadlerEtAl15}
{\sc T.~Isaac, N.~Petra, G.~Stadler, and O.~Ghattas}, {\em Scalable and
  efficient algorithms for the propagation of uncertainty from data through
  inference to prediction for large-scale problems, with application to flow of
  the {A}ntarctic ice sheet}, Journal of Computational Physics, 296 (2015),
  pp.~348--368.

\bibitem{KazeevKhammashNipSchwab14}
{\sc V.~Kazeev, M.~Khammash, M.~Nip, and C.~Schwab}, {\em Direct solution of
  the chemical master equation using quantized tensor trains}, PLoS
  computational biology, 10 (2014).

\bibitem{KazeevReichmannSchwab13}
{\sc V.~Kazeev, O.~Reichmann, and C.~Schwab}, {\em Low-rank tensor structure of
  linear diffusion operators in the tt and qtt formats}, Linear Algebra and its
  Applications, 438 (2013), pp.~4204--4221.

\bibitem{Khoromskij11}
{\sc B.~N. Khoromskij}, {\em O(d log {N})-quantics approximation of {N}-d
  tensors in high-dimensional numerical modeling}, Constructive Approximation,
  34 (2011), pp.~257--280.

\bibitem{KoldaBader09}
{\sc T.~G. Kolda and B.~W. Bader}, {\em Tensor decompositions and
  applications}, SIAM Review, 51 (2009), pp.~455--500.

\bibitem{KoldaMayo11}
{\sc T.~G. Kolda and J.~R. Mayo}, {\em Shifted power method for computing
  tensor eigenpairs}, SIAM Journal on Matrix Analysis and Applications, 32
  (2011), pp.~1095--1124.

\bibitem{LinLuYing11}
{\sc L.~Lin, J.~Lu, and L.~Ying}, {\em Fast construction of hierarchical matrix
  representation from matrix--vector multiplication}, Journal of Computational
  Physics, 230 (2011), pp.~4071--4087.

\bibitem{LoggMardalGarth12}
{\sc A.~Logg, K.-A. Mardal, and G.~Wells}, {\em Automated Solution of
  Differential Equations by the Finite Element Method: {T}he {FEniCS} book},
  vol.~84, Springer Science \& Business Media, 2012.

\bibitem{Maddison19}
{\sc J.~R. Maddison, D.~N. Goldberg, and B.~D. Goddard}, {\em Automated
  calculation of higher order partial differential equation constrained
  derivative information}, SIAM Journal on Scientific Computing, 41 (2019),
  pp.~C417--C445.

\bibitem{Naumann12}
{\sc U.~Naumann}, {\em The Art of Differentiating Computer Programs: An
  Introduction to Algorithmic Differentiation}, Society for Industrial and
  Applied Mathematics, 2012.

\bibitem{Oseledets09}
{\sc I.~V. Oseledets}, {\em Approximation of matrices with logarithmic number
  of parameters}, in Doklady Mathematics, vol.~80, Springer, 2009,
  pp.~653--654.

\bibitem{Oseledets10}
\leavevmode\vrule height 2pt depth -1.6pt width 23pt, {\em Approximation of
  $2^{d}\times 2^{d}$ matrices using tensor decomposition}, SIAM Journal on
  Matrix Analysis and Applications, 31 (2010), pp.~2130--2145.

\bibitem{Oseledets11}
\leavevmode\vrule height 2pt depth -1.6pt width 23pt, {\em Tensor-train
  decomposition}, SIAM Journal on Scientific Computing, 33 (2011),
  pp.~2295--2317.

\bibitem{OseledetsTyrtyshnikov10}
{\sc I.~V. Oseledets and E.~Tyrtyshnikov}, {\em {TT}-cross approximation for
  multidimensional arrays}, Linear Algebra and its Applications, 432 (2010),
  pp.~70--88.

\bibitem{OseledetsTyrtyshnikovZamarashkin11}
{\sc I.~V. Oseledets, E.~Tyrtyshnikov, and N.~Zamarashkin}, {\em Tensor-train
  ranks for matrices and their inverses}, Computational Methods in Applied
  Mathematics, 11 (2011), pp.~394--403.

\bibitem{PerezEtAl06}
{\sc D.~Perez-Garcia, F.~Verstraete, M.~Wolf, and J.~Cirac}, {\em Matrix
  product state representations}, Quantum Information \& Computation, 7 (2007),
  pp.~401--430.

\bibitem{PetraMartinStadlerEtAl14}
{\sc N.~Petra, J.~Martin, G.~Stadler, and O.~Ghattas}, {\em A computational
  framework for infinite-dimensional {B}ayesian inverse problems: {P}art {II}.
  {S}tochastic {N}ewton {MCMC} with application to ice sheet inverse problems},
  SIAM Journal on Scientific Computing, 36 (2014), pp.~A1525--A1555.

\bibitem{SavostyanovOseledets11}
{\sc D.~Savostyanov and I.~Oseledets}, {\em Fast adaptive interpolation of
  multi-dimensional arrays in tensor train format}, in The 2011 International
  Workshop on Multidimensional (nD) Systems, IEEE, 2011, pp.~1--8.

\bibitem{SpantiniSolonenCuiEtAl15}
{\sc A.~Spantini, A.~Solonen, T.~Cui, J.~Martin, L.~Tenorio, and Y.~Marzouk},
  {\em Optimal low-rank approximations of {B}ayesian linear inverse problems},
  SIAM Journal on Scientific Computing, 37 (2015), pp.~A2451--A2487.

\bibitem{YangKrompassTresp17}
{\sc Y.~Yang, D.~Krompass, and V.~Tresp}, {\em Tensor-train recurrent neural
  networks for video classification}, in Proceedings of the 34th International
  Conference on Machine Learning-Volume 70, JMLR. org, 2017, pp.~3891--3900.

\bibitem{ZhuLiFomelEtAl16}
{\sc H.~Zhu, S.~Li, S.~Fomel, G.~Stadler, and O.~Ghattas}, {\em A {B}ayesian
  approach to estimate uncertainty for full waveform inversion with a priori
  information from depth migration}, Geophysics, 81 (2016), pp.~R307--R323.

\end{thebibliography}
	
\end{document}